\documentclass[12pt]{article}

\usepackage{subcaption}

\usepackage{float}

\usepackage{multirow}

\usepackage{tikz}

\usetikzlibrary{calc,backgrounds,arrows,matrix}

\usepackage{enumerate}

\usepackage{color}
\definecolor{lightblue}{rgb}{0,0.2,0.5}

\usepackage[colorlinks=true, urlcolor=blue,linkcolor=blue, citecolor=lightblue]{hyperref}

\usepackage[round,sort,comma,numbers]{natbib}

\bibpunct{\textcolor{lightblue}{(}}{\textcolor{lightblue}{)}}{,}{a}{}{;}

\usepackage{amssymb,amsmath}

\usepackage{graphicx}

\usepackage{graphics}
\usepackage{color}

\usepackage{tikz}
\usetikzlibrary{trees}

\DeclareMathAlphabet{\eufrak}{U}{}{}{}
\SetMathAlphabet\eufrak{normal}{U}{euf}{m}{n}
\SetMathAlphabet\eufrak{bold}{U}{euf}{b}{n}

\allowdisplaybreaks

 \def\qu{{\mathord{\mathbb Z}}}

 \def\T{{\mathrm{{\rm T}}}}

 \def\sZZ{{\rm Z\kern-.45em{}Z}}

 \def\sQQ{{\kern 0.27em \vrule height1.45ex width0.03em depth0em
           \kern-0.30em \rm Q}}
 \def\qu{{\mathchoice
         {\sQQ}
         {\sQQ}
   {\kern 0.225em \vrule height1.05ex width0.025em depth0em \kern-0.25em \rm Q}
   {\kern 0.180em \vrule height0.78ex width0.020em depth0em \kern-0.20em \rm Q}
         }}
 \def\sGG{{\kern 0.27em \vrule height1.45ex width0.03em depth0em
           \kern-0.30em \rm G}}
 \def\gg{{\mathchoice
         {\sGG}
         {\sGG}
   {\kern 0.225em \vrule height1.05ex width0.025em depth0em \kern-0.25em \rm G}
   {\kern 0.180em \vrule height0.78ex width0.020em depth0em \kern-0.20em \rm G}
         }}

 \newtheorem{prop}{Proposition}[section]
 
 \newtheorem{definition}[prop]{Definition}
 \newtheorem{corollary}[prop]{Corollary}
 \newtheorem{theorem}[prop]{Theorem}
 \newtheorem{remark}[prop]{Remark}
 \newtheorem{example}[prop]{Example}

\numberwithin{equation}{section}

 \def\P{{\mathord{\mathbb P}}}

 \newcounter{hyp}
\usepackage{aeguill}

\newenvironment{Proof}{\removelastskip\par\medskip \noindent{\em Proof.} \rm}{\penalty-20\null\hfill$\square$\par\medbreak}

\def\bprf{\begin{Proof}}
\def\nprf{\end{Proof}}
\def\bdes{\begin{description}}
\def\ndes{\end{description}}

\newtheorem{thm}{Theorem}[section]

\newcommand{\sgn}{\mathrm{sign }\:}

\def\bdef{\begin{defn}}
\def\ndef{\end{defn}}
\def\bthm{\begin{thm}}
\def\nthm{\end{thm}}
\def\bprop{\begin{prop}}
\def\nprop{\end{prop}}
\def\brmk{\begin{remark}}
\def\nrmk{\end{remark}}
\def\bexa{\begin{exa}}
\def\nexa{\end{exa}}
\def\blem{\begin{lem}}
\def\nlem{\end{lem}}
\def\bcor{\begin{cor}}
\def\ncor{\end{cor}}
\def\bexe{\begin{exe}}
\def\nexe{\end{exe}}

\def\R{\eufrak{R}}

\newcommand{\real}{\mathbb{R}}

\def\og{\leavevmode\raise.3ex
     \hbox{$\scriptscriptstyle\langle\!\langle$~}}
\def\fg{\leavevmode\raise.3ex
     \hbox{~$\!\scriptscriptstyle\,\rangle\!\rangle$}~}

\title{\Huge
  {Existence} and probabilistic representation of the solutions of semilinear
  {parabolic} PDEs with fractional Laplacians
}

\author{
  Guillaume Penent\footnote{\href{mailto:PENE0001@e.ntu.edu.sg}{pene0001@e.ntu.edu.sg}}
  \qquad
      Nicolas Privault\footnote{
\href{mailto:nprivault@ntu.edu.sg}{nprivault@ntu.edu.sg}
}
  \\
\small
Division of Mathematical Sciences
\\
\small
School of Physical and Mathematical Sciences
\\
\small
Nanyang Technological University
\\
\small
21 Nanyang Link, Singapore 637371
}

\allowdisplaybreaks

\usepackage{empheq}

\usepackage{bbm}

\makeatletter
\newcommand*\rel@kern[1]{\kern#1\dimexpr\macc@kerna}
\newcommand*\widebar[1]{
  \begingroup
  \def\mathaccent##1##2{
    \rel@kern{0.8}
    \overline{\rel@kern{-0.8}\macc@nucleus\rel@kern{0.2}}
    \rel@kern{-0.2}
  }
  \macc@depth\@ne
  \let\math@bgroup\@empty \let\math@egroup\macc@set@skewchar
  \mathsurround\z@ \frozen@everymath{\mathgroup\macc@group\relax}
  \macc@set@skewchar\relax
  \let\mathaccentV\macc@nested@a
  \macc@nested@a\relax111{#1}
  \endgroup
}
\makeatother

\let\oldcitet=\citet
\let\oldcitep=\citep
\renewcommand{\cite}[1]{\textcolor[rgb]{0,0,1}{\oldcitet{#1}}}
\renewcommand{\citet}[1]{\textcolor[rgb]{0,0,1}{\oldcitet{#1}}}
\renewcommand{\citep}[1]{\textcolor[rgb]{0,0,1}{\oldcitep{#1}}}

\begin{document}

\maketitle

\baselineskip0.6cm

\vspace{-0.6cm}

\begin{abstract}
{
  We obtain existence results for the solution $u$ of nonlocal semilinear
  parabolic PDEs on $\real^d$ with polynomial nonlinearities in $(u, \nabla u)$, using} a tree-based probabilistic representation. This probabilistic representation applies to the solution of the equation itself, as well as to its partial derivatives by associating one of $d$ marks to the initial tree branch. {Partial derivatives} are dealt with by integration by parts and subordination of Brownian motion. Numerical illustrations are provided {in examples for the fractional Laplacian in dimension up to $10$, and for the fractional Burgers equation in dimension two}.
\end{abstract}



\baselineskip0.7cm

\parskip-0.1cm

\section{Introduction}
Nonlocal partial differential operators such as the fractional
Laplacian are useful in the modeling of anomalous diffusion phenomena
driven in particular by stable L\'evy processes,
and
they have found applications in multiple fields of
engineering, physics and finance.
The numerical solution of elliptic boundary value problems
involving fractional Laplacians have been studied
{
by means of finite differences in 
the one-dimensional case in e.g. \cite{oberman} in the parabolic
case, and in \cite{acosta2}, and \cite{acosta2021}
 in the elliptic case.
}

\medskip

\indent
Probabilistic approaches relying on the Feynman-Kac formula
represent alternatives to finite differences for the numerical
solution of parabolic partial differential equations. 
The use of stochastic diffusion branching mechanisms
for the representation of solutions of partial differential
equations has been introduced by \cite{skorohodbranching},
and this construction
has been extended in \cite{inw} to branching Markov processes.
In \cite{N-S}, branching Markov processes have been applied
to the blowup of solutions of a wide class of parabolic PDEs
using their Duhamel integral formulations
and the Markov property of the branching process
at its first branching time.
The branching mechanism has also been applied in \cite{hpmckean}
to the KPP equation, and
to the blow-up of solutions of \cite{fujita} equations
of the form $\partial u(t,x)/\partial t = \Delta u(t,x) + c u^\beta (t,x)$
in \cite{lm}, 
see also \cite{chakraborty} 
for the existence of solutions
of parabolic PDEs with 
power series nonlinearities.
Related arguments have also been applied to Fourier-transformed PDEs
in order to treat the Navier-Stokes equation by the
use of stochastic cascades in \cite{sznitman},
see also \cite{blomker}
for the representation of Fourier modes for the
solution of class of semilinear parabolic PDEs.

\medskip

\indent
 This branching argument has been recently extended in
 \cite{labordere} to the treatment polynomial non-linearities in
 gradient terms.
 {
     For this, branches associated to gradient terms
 are specified using marks,
 and are subject to Malliavin integrations by parts.} 
   This approach applies in principle to continuous It\^o diffusion generators,
provided that the corresponding Malliavin weight can be successfully estimated.
{In the absence of gradient nonlinearities,
  the tree-based approach has been recently implemented
  for nonlocal semilinear PDEs in \cite{belak}.
  } 

\medskip

\indent
 {In this paper, 
 we obtain existence results for the solution of 
 nonlocal semilinear PDEs
 by extending the above arguments from 
 the standard Laplacian $\Delta$ to}  
pseudo-differential operators of the form $-\eta(-\Delta /2)$, where
$\eta$ is a Bernstein function such that $\eta(0^+) = 0$.
 
\medskip

 Precisely, given a horizon time $T>0$, we consider the semilinear PDE
 given as
 \begin{equation}
   \label{eq:1}
\begin{cases}
  \displaystyle
  \frac{\partial u}{\partial t} (t,x)
  - \eta ( - \Delta / 2 ) u(t,x) +
  f\left(t, x,u(t,x) , \frac{\partial u}{\partial x_1}(t,x), \ldots ,
  \frac{\partial u}{\partial x_m}(t,x) \right)
  = 0,
  \medskip
  \\
u(T,x) = \phi(x), \qquad x=(x_1,\ldots , x_d) \in \real^d,
\end{cases}
\end{equation}
 where 
 $f(t,x,y,z_1,\ldots , z_m)$ is a polynomial nonlinearity given by
$$
f(t, x,y , z_1,\ldots , z_m) =
\sum\limits_{l= (l_0,\ldots , l_m) \in {\cal L}_m} c_l(t,x) y^{l_0}
z^{l_1}_1 \cdots z^{l_m}_m,
$$
$t\in \real_+$, $x,y,z_1,\ldots , z_m\in \real$,
for some $m\in \{0,\ldots ,d \}$,
 where ${\cal L}_m$ is a finite subset of $\mathbb{N}^{m+1}$
 and
$c_l(t,x)$ are measurable functions of $(t,x)\in [0,T]\times \real^d$,
$l= (l_0,\ldots , l_m)\in {\cal L}_m$.
{ In the sequel, we let $\Vert x\Vert := \sqrt{x_1^2+\cdots + x_d^2}$,
 $x=(x_1,\ldots , x_d) \in \real^d$. 
}
\\
 \medskip\noindent
{\bf Assumption (\hypertarget{BGJhyp}{$A$}):} 
{\em
 We assume that the coefficients $c_l (t,x)$ are uniformly bounded,
 i.e.
 \begin{equation}
   \label{cc1}
   |c_l|_{\infty} := \sup_{t \in [0,T], x \in \mathbb{R}^d} |c_l(t,x)| < \infty,
 \qquad l = (l_0,\ldots ,l_m) \in {\cal L}_m,
 \end{equation}
  and that the terminal condition $\phi$ is Lipschitz, i.e.
  \begin{equation}
    \label{cc2}
    |\phi(x)-\phi(y)| \leq L \Vert x-y\Vert,
    \qquad
     x,y \in \real^d,
  \end{equation}
  for some $L>0$, and bounded on $\real^d$.
}
\medskip
 
 \noindent 
  In the sequel, we will say that a function $u(t,x)$ is an 
 {integral} solution
  if $u(t,x)$  satisfies the Duhamel formulation of \eqref{eq:1}, i.e.
 \begin{eqnarray}
   \label{fjkdslf}
    u (t,x) & = & \int_{\real^d} \varphi (T-t,y-x)\phi(y)dy
   \\
   \nonumber
   & &
   +
   \sum_{l = (l_0,\ldots ,l_m) \in {\cal L}_m }
   \int_t^T \int_{\real^d} \varphi (s-t,y-x) c_l(s,y) u^{l_0}(s,y) \prod_{j=1}^m \left(
   \frac{\partial u}{\partial y_j} (s,y)\right)^{l_j} dy ds,
\end{eqnarray}
 $(t,x) \in [0,T]\times \real^d$.
 Note that the above setting includes the
 case of the standard
 fractional Laplacian $\Delta_\alpha  = - \eta ( - \Delta / 2 )$
 by choosing the Laplace exponent $\eta(\lambda) =  (2\lambda)^{\alpha/2}$.

 \medskip 
 
In particular, in Theorem~\ref{t1} we provide probabilistic representations for
the solutions of a wide class of semilinear parabolic PDEs of the form
\begin{equation}
  \label{fjkldsf}
 \frac{\partial u}{\partial t}(t,x) - \eta(-\Delta /2)u(t,x)
 +
  f\left(t, x,u(t,x) , \frac{\partial u}{\partial x_1}(t,x), \ldots ,
  \frac{\partial u}{\partial x_m}(t,x) \right)
 = 0, ~~ u(T,\cdot ) = \phi (\cdot ),
\end{equation}
 $(t,x) \in [0,T] \times \mathbb{R}^d$,
 with polynomial non-linearity $f$ in the solution $u$
 and its {
   partial derivatives $\partial u / \partial x_i$,
 $i=1,\ldots , m$, } 
 and $\eta$ is a Bernstein function that satisfies $\eta(0^+)=0$.

 \medskip
 
 {
 The probabilistic representations of 
 Theorem~\ref{t1} uses a functional
 of a random branching process driven
 by a subordinated L\'evy process
 $(Z_t)_{t\in \real_+} := (B_{S_t})_{t\in \real_+}$, where
 $(B_t)_{t\in \real_+}$ is a standard $d$-dimensional Brownian motion
 and $(S_t)_{t\in \real_+}$
 is a L\'evy subordinator with Laplace exponent $\eta$ 
 such that 
$$
\mathbb{E}\big[e^{-\lambda S_t}\big] = e^{-t \eta( \lambda )}, \qquad
\lambda , t \geq 0, 
$$
see e.g. Theorem~1.3.23 and pages~55-56 in \cite{applebk2}.
 Then, by Proposition~1.3.27 in \cite{applebk2},
   $(Z_t)_{t\in \real_+}$ has L\'evy symbol 
$\psi_Z (\xi) = -\eta( \Vert\xi \Vert^2/2 ) $ such that 
  $\mathbb{E}\big[e^{ i \xi Z_t}\big] = e^{t \psi_Z ( \xi )}$,
 $\xi\in \real^d$, $t \geq 0$, 
 and, by Theorem~3.3.3 therein, 
 the infinitesimal generator 
 of $(Z_t)_{t\in \real_+}$ is
the pseudo-differential operator $-\eta(-\Delta / 2 )$. 
}

 \medskip

 In the case of stable processes we have
 $\eta(\lambda) := ( 2 \lambda)^{\alpha / 2}$, and
 $-\eta(-\Delta / 2 )$ becomes the fractional Laplacian
$$
 \Delta_\alpha u = - ( - \Delta )^{\alpha / 2} u
 = \frac{2^\alpha\Gamma(d/2+\alpha/2)}{\pi^{d/2}|\Gamma(-\alpha/2)|}
 \lim_{r\rightarrow 0^+}
 \int_{\mathbb{R}^d \backslash B(x,r)} \frac{u( \cdot +z)-u(z)}{|z|^{d+\alpha}}dz,
 $$
 for $\alpha \in (0,2]$, where $\Gamma(p) : = \displaystyle
   \int_0^\infty e^{-\lambda x} \lambda^{p-1} d\lambda$, $p>0$,
   is the gamma function, see e.g. \cite{tendef}.
 
\medskip 

 For each $i=0,1,\ldots ,d$
 we construct a sufficiently integrable functional
 $\mathcal{H}_\phi (\mathcal{T}_{t,x,i})$
 of a random tree $\mathcal{T}_{t,x,i}$ such that we have
 the representations
 $$
 u (t,x) :=  \mathbb{E}\big[ \mathcal{H}_\phi (\mathcal{T}_{t,x,0})\big],
 \quad (t,x) \in [0,T]\times \real^d,
 $$
 and
$$
\frac{\partial u}{\partial x_i} (t,x) :=  \mathbb{E}\big[ \mathcal{H}_\phi (\mathcal{T}_{t,x,i})\big],
\quad (t,x) \in [0,T]\times \real^d,
\quad i=1,\ldots ,d,
$$
see Theorem~\ref{t1}.
{Dealing with gradient terms in the proof of Theorem~\ref{t1}
  requires to perform an integration by parts,
  which is made possible using the Gaussian density of $B_t$
 in the subordination
 $Z_t := B_{S_t}$,} as done in
\cite{kawai3} in the case of stable processes
with $\eta(\lambda) := ( 2 \lambda)^{\alpha / 2}$. 

   \medskip

 As a consequence of Theorem~\ref{t1}, in Proposition~\ref{fkdslf1}
 we show that the probabilistic representation of
 Theorem~\ref{t1} can be used to recover
 the classical result of \cite{fujita} on the blow-up of semilinear PDEs.  

   \medskip

 While the branching tree mechanism is quite general and
can be applied to a wide range of differential equations via formal
calculations, proving the existence of solutions requires to
show the integrability of
{ functional $\mathcal{H}_\phi (\mathcal{T}_{t,x,i})$
  representing the PDE solution $u(t,x)$ and its partial derivatives.
 We deal with this integrability 
 using existence results for the}
 solutions of Volterra integral equations, 
 instead of using ODEs as in e.g. \cite{labordere2} and \cite{labordere}.

 \medskip

 Theorem~\ref{t2}, we show 
 that the integrability required for the probabilistic representation
 Theorem~\ref{t1} is satisfied provided that
 $\lambda \rightarrow 1 / ( \sqrt{\lambda} \eta (\lambda) )$ is
 integrable at $+\infty$.
 In comparison with recent work in the diffusion case, see \cite{labordere},
 {
 our integrability condition \eqref{eq:IC}-\eqref{eq:IC2} in Theorem~\ref{t1} 
 is sharper because it only involves mark indexes of
 partial derivatives appearing in the main PDE.
 In addition, we provide a detailed justification for
 the commutation relation \eqref{fskl2}
 instead of stating it as an assumption as in \cite{labordere},
 see Assumption 3.2 therein.
 }
 
 \medskip
 
{
 As a direct consequence of Theorems~\ref{t1} and
 \ref{t2}, we obtain the following result on
   local-in-time existence of solutions.
  \begin{theorem}
 Under Assumption~(\hyperlink{BGJhyp}{$A$}), suppose that 
 $$
 \int_{\lambda_0}^\infty
 \frac{1}{\sqrt{\lambda} \eta (\lambda)} d\lambda < \infty
 $$
 for some $\lambda_0 > 0$.
 Then, there exists
 a small enough $T>0$ such that
 the PDE \eqref{eq:1}
 admits an integral solution on $[0,T]$ in the sense of \eqref{fjkdslf}.
 \end{theorem}
   Related local and global-in-time existence results have been obtained
 for generalized fractional Laplacians by deterministic arguments
 under more technical conditions 
 in e.g. \cite{ishige2014} and more recently in \cite{ishige2021}
 for power nonlinearities
 of sufficiently low orders.
 In the particular case of the $\alpha$-fractional Laplacian 
 where $\eta(\lambda) := (2\lambda)^{\alpha / 2}$ with $\alpha \in (1,2)$, 
 we obtain the following corollary. 
\begin{corollary}
 Taking $\eta(\lambda) := (2\lambda)^{\alpha / 2}$
  with $\alpha \in (1,2)$, 
  under Assumption~(\hyperlink{BGJhyp}{$A$})
 there exists
 a small enough $T>0$ such that
 the PDE \eqref{eq:1} with $\alpha$-fractional Laplacian 
  admits an integral solution on $[0,T]$ in the sense of \eqref{fjkdslf}. 
\end{corollary} 
\noindent
  In the case of the fractional Laplacian, 
 Proposition~\ref{t3}
 provides quantitative estimates
 on the horizon time $T$, 
 ensuring existence of solutions on $[0,T]$ by Theorem~\ref{t1}. 
}

\medskip

{
  We also provide a Monte Carlo implementation of our algorithm
  for the numerical solutions
  of nonlinear fractional PDEs with and without gradient term
   in dimension up to 10, 
  and of a fractional Burgers equation.
 The tree-based Monte Carlo method 
 avoids the curse of dimensionality, 
 whereas the application of deterministic numerical methods
 is notoriously difficult
 including in the fractional case, see, e.g., \cite{bonito}. 
 } 

\medskip

 The paper is organized as follows.
In Section~\ref{s2} we present
 the description of the branching mechanism in Section~\ref{s2}.
In Section~\ref{s3} we state our main result
 Theorem~\ref{t1}
 which gives the probabilistic representation of the solution and
 its partial derivatives.
 In Section~\ref{s4}
 we give give a sufficient condition on
 the Bernstein function $\eta$
 that ensures the integrability needed for the
 the probabilistic representation of Theorem~\ref{t1} to hold.
 In Section~\ref{s5},
 we present some numerical simulations to illustrate the method on specific examples.

 \subsubsection*{Bernstein functions and subordinators}
Let $\eta :(0,\infty) \rightarrow [0,\infty) $
  denote a Bernstein function,
i.e. $\eta$ is a $C^\infty$ function
whose $nth$ derivative satisfies $(-1)^n \eta^{(n)} \leq 0$ for all $n \geq 1$,
and $\lim_{z\searrow 0} \eta(z) = 0$,
see Theorem~1.3.23 in \cite{applebk2}.
We consider a subordinator $(S_t)_{t\in \real_+}$,
i.e. a $\real_+$-valued non-decreasing L\'evy process,
with Laplace exponent $\eta$, which admits the representation
 \begin{equation}
 \label{djkldsf}
\eta (\lambda) =  b\lambda +\int_0^\infty(1-e^{-\lambda y})\nu(dy),
\end{equation}
where $b \geq 0$ and the L\'evy measure $\nu$ satisfies
$$
\int_0^\infty (y \wedge 1) \nu(dy) < \infty,
$$
see Theorem~1.3.15 in \cite{applebk2}. 
Using the identity 
\begin{equation}
\label{id}
  x^{-p}= \frac{1}{\Gamma(p)}\int_0^\infty e^{-\lambda x} \lambda^{p-1} d\lambda
   \qquad x>0, 
  \end{equation}
 the negative moments of $S_t$ are given by 
\begin{equation} 
\label{id2}
  \mathbb{E}\big[ S_t^{-p} \big]
  = \frac{1}{\Gamma(p)}\int_0^\infty e^{-t \eta(\lambda)} \lambda^{p-1} d\lambda,
  \qquad p >0. 
\end{equation} 

\noindent 
When $(S_t)_{t\in \real_+}$ is an $\alpha/2$-stable subordinator with Laplace exponent $\eta(\lambda) = (2\lambda)^{\alpha/2}$,
the subordinated process
$Z_t=B_{S_t}$ becomes an $\alpha$-stable process with generator $\Delta_\alpha $. In that case, we have $b=0$ in \eqref{djkldsf}, the L\'evy measure $\nu$ of the subordinator $(S_t)_{t\in \real_+}$ is given by
$$
\nu(dx) = \alpha \frac{2^{\alpha/2 -1 }}{\Gamma(1-\alpha/2)} \frac{dx}{x^{1+\alpha/2}},
$$
 and its L\'evy symbol $\psi_S$ satisfies
\begin{eqnarray}
  \nonumber
  \psi_S (\xi) &= & \frac{2^{\alpha/2}\alpha}{2\Gamma(1-\alpha/2)}\int_0^\infty (e^{i\xi y}-1)\frac{dy}{y^{1+\alpha/2}}
  \\
  \nonumber
  &=& \frac{\alpha (2|\xi |)^{\alpha/2}}{2\Gamma(1-\alpha/2)}\Gamma(-\alpha/2)e^{i
  \alpha \arg(-i\xi ) / 2}\\
\label{ls}
  &= & -\cos\left( \frac{\pi \alpha }{4}\right)
(2|\xi |)^{\alpha/2}
\left( 1-i \ \! \sgn(\xi )\tan\left(
\frac{\pi \alpha}{4}\right)\right), \qquad \xi \in \real,
\end{eqnarray}
where we used the identity
$$
\int_0^\infty (e^{wy}-1) y^{-1-\alpha/2} dy = \Gamma(-\alpha/2) |w|^{\alpha/2} e^{i \alpha \arg(-w) / 2 }
$$
which is valid for $\alpha \in (0,2)$ and any $w \in \mathbb{C}^*$ with $\R (w) \leq 0$, see Relation~(14.18) page~84 of \cite{sato}.
 In this case, the negative moments of $S_t$ are given by 
\begin{eqnarray}
  \nonumber
  \mathbb{E}\big[ S_t^{-p} \big]
  &=& \frac{1}{\Gamma(p)}\int_0^\infty e^{-t  (2\lambda)^{\alpha/2}} \lambda^{p-1} d\lambda 
\\
\nonumber
  &=& \frac{1}{t^{2p/\alpha}}\frac{2^{1-p}}{\alpha\Gamma(p)}\int_0^\infty u^{-1+2p/\alpha } e^{-u}du\\
\label{id3}
&=& \frac{2^{1-p} \Gamma(2p/\alpha)}{\alpha t^{2p/\alpha} \Gamma(p)}, \qquad p>0. 
\end{eqnarray}
\section{Random trees with marked branches}
\label{s2}
In the sequel, we will provide a probabilistic representation for the solution of \eqref{eq:1}, using a branching mechanism
such that the solution of \eqref{eq:1}
will be given by the expectation of a multiplicative functional defined
on a random tree structure.

\medskip

Let $\rho: \mathbb{R}^+ \rightarrow (0,\infty )$ be a
probability density function on $\real_+$,
and consider a probability mass function
$(q_{l_0,\ldots ,l_m})_{(l_0,\ldots ,l_m)\in {\cal L}_m}$ on ${\cal L}_m$
 with $q_{l_0,\ldots ,l_m} > 0$,
$(l_0,\ldots ,l_m)\in {\cal L}_m$, and
$\sum_{ (l_0,\ldots ,l_m) \in {\cal L}_m} |l| q_{l_0,\ldots ,l_m} < \infty$,
where $|l| = l_0+\cdots +l_m$.
 In addition, we consider 
\begin{itemize}
\item an i.i.d. family $(\tau^{i,j})_{i,j\geq 1}$ of random variables
 with distribution $\rho (t)dt$ on $\real_+$,
\item an i.i.d. family $(I^{i,j})_{i,j\geq 1}$ of discrete
  random variables, with
  \begin{equation}
    \label{fjdks}
  \mathbb{P}\big( I^{i,j}=(l_0,\ldots ,l_m) \big) = q_{l_0,\ldots ,l_m} >0,
  \qquad (l_0,\ldots , l_m)\in {\cal L}_m,
\end{equation}
\item
  an independent family $(Z^{i,j})_{i,j\geq 1}$
  of subordinated L\'evy processes
  constructed as
  $$
  Z^{i,j}_t := B^{i,j}_{S^{i,j}_t},
  \qquad
  t\in \real_+, \quad i,j \geq 1,
  $$
  where
  $(B^{i,j})_{i,j\geq 1}$
  and
  $(S^{i,j})_{i,j\geq 1}$
  are independent standard Brownian motions
  and independent subordinators with
  Laplace exponent $\eta$.
\end{itemize}
 The sequences $(\tau^{i,j})_{i,j\geq 1}$, $(I^{i,j})_{i,j\geq 1}$ and
 $(Z^{i,j})_{i,j\geq 1}$ are assumed to be mutually independent.
\subsubsection*{Marked branching process}
We consider a marked branching process starting
from a particle at the position $x\in \real^d$, with label $\widebar{1}=(1)$
and mark $i \in \{0,1,\ldots ,d \}$ at time $t\in [0,T]$,
which evolves according to 
the process $X_{s,x}^{\widebar{1}} = x + Z_{s-t}^{1,1}$, $s \in [t,t+\tau^{1,1}]$.

\medskip

If $\tau^{1,1}<T-t$, the process branches at time $t+\tau^{1,1}$
into new independent copies of
$(Z_t)_{t \in \real_+}$, each of them started at
 the position $X_{t+\tau^{1,1}}$ at time $t+\tau^{1,1}$.
Based on the values of $I^{1,1} =(l_0,\ldots , l_m)\in {\cal L}_m$,
a family of $|l|:=l_0+\cdots +l_m$ of new branches
carrying respectively the marks $i=0,\ldots ,d$
are created with the probability $q_{l_0,\ldots ,l_m}$,
where
\begin{itemize}
\item the first $l_0$ branches
  carry the mark $0$ and
  are indexed by $(1,1),(1,2),\ldots ,(1,l_0)$,
\item the next $l_1$ branches
  carry the mark $1$ and
  are indexed by $(1,l_0+1),\ldots ,(1,l_0 + l_1)$, and so on.
\end{itemize}
Each new particle then follows independently
the same mechanism as the first one, and
every branch stops when it reaches the horizon time $T$.
Particles at the generation $n\geq 1$ are assigned a label of the form
 $\widebar{k} = (1,k_2,\ldots ,k_n) \in \mathbb{N}^n$,
and their parent is labeled $\widebar{k}- := (1,k_2,\ldots ,k_{n-1})$.
The particle labeled $\widebar{k}$ is born at time $T_{\widebar{k}-}$
and its lifetime $\tau^{n,\pi_n(\widebar{k})}$ is the element of index
$\pi_n(\widebar{k})$ in the i.i.d. sequence
$(\tau^{n,j})_{j\geq 1}$,
defining a random injection
$$\pi_n:\mathbb{N}^n \to \mathbb{N},
\qquad n\geq 1.
$$
The random evolution of particle $\widebar{k}$
is given by
\begin{equation}
  \label{fdjskldf}
X_{s,x}^{\widebar{k}} := X^{\widebar{k}-}_{T_{\widebar{k}-},x}+Z_{s-T_{\widebar{k}-}}^{n,\pi_n(\widebar{k})},
\qquad s\in [T_{\widebar{k}-},T_{\widebar{k}}],
\end{equation}
  where $T_{\widebar{k}} := \T_{\widebar{k}-} + \tau^{n,\pi_n(\widebar{k})}$.

\medskip

 If $T_{\widebar{k}} := \T_{\widebar{k}-} + \tau^{n,\pi_n(\widebar{k})} < T$,
 we draw a sample
 $I_{\widebar{k}}: = I^{n,\pi_n(\widebar{k})} = (l_0,\ldots ,l_m)$
 of $I^{n,\pi_n(\widebar{k})}$ with distribution \eqref{fjdks},
 and the particle $\widebar{k}$ branches into
 $|I^{n,\pi_n(\widebar{k})}|=l_0+\cdots +l_m$ offsprings at generation $(n+1)$,
 which are indexed by $(1,\ldots ,k_n,i)$, $i=1,\ldots ,|I^{n,\pi_n(\widebar{k})}|$.
 The particles whose index ends with an integer between $1$ and $l_0$
 will carry the mark $0$, and those with index ending with an integer between
 $l_0+\cdots +l_{i-1} +1$ and $l_0+\cdots + l_i$ will carry a mark
 $i\in \{1,\ldots ,d\}$.
 Finally, the mark of particle $\widebar{k}$ will be denoted by
 $\theta_{\widebar{k}} \in \{0,1,\ldots ,d \}$.
Note that the indexes are only be used to distinguish the particles in the
branching process, and they are distinct from the marks.

 \medskip

 The set of particles dying before time $T$ is denoted by $\mathcal{K}^{\circ}$,
whereas those dying after $T$ form a set denoted by $\mathcal{K}^{\partial}$.
\begin{definition}
 When started at time $t\in [0,T]$ from
 a position $x\in \real^d$ and a mark
 $i \in \{0,1,\ldots ,d\}$ on its first branch,
 the above construction yields
 a marked branching process called a random marked tree, and
 denoted by $\mathcal{T}_{t,x,i}$.
\end{definition}
 The tree $\mathcal{T}_{t,x,0}$
 will be used for the stochastic representation of the solution
 $u$ of the PDE \eqref{eq:1}, while the trees $\mathcal{T}_{t,x,i}$
 will be used for the stochastic representation of the
 partial derivatives $\partial u/\partial x_i$,
 $i=1,\ldots , m$.
 The next table summarizes the notation introduced so far.

 \medskip
 \bigskip

\begin{center}
 \begin{tabular}{||l | c||}
 \hline
 Object & Notation \\ [0.5ex]
 \hline\hline
 Initial time  & $t$ \\
 \hline
 Initial position & $x$ \\
 \hline
 Tree rooted at $(t,x)$ with initial mark $i$ &  $\mathcal{T}_{t,x,i}$ \\  \hline
 Particle (or label) of generation $n\geq 1$ & $\stackrel{}{\widebar{k}=(1,k_2,\ldots ,k_n)}$\\
 \hline
 First branching time & $T_{\widebar{1}}$\\
 \hline
 Lifespan of a particle & $\tau_{\widebar{k}} = T_{\widebar{k}} - T_{\widebar{k}-}$ \\
 \hline
Birth time of a particle $\bar{k}$ & $T_{\widebar{k}-}$ \\
\hline
Death time of a particle $\bar{k}$ & $T_{\widebar{k}}$ \\
\hline
Position at birth & $X^{\widebar{k}}_{T_{\widebar{k}-},x}$\\
\hline
Position at death & $X^{\widebar{k}}_{T_{\widebar{k}},x}$ \\
\hline
Mark & $\theta_{\widebar{k}}$
\\ 
\hline
\end{tabular}
\end{center}

\bigskip
\noindent
To represent the structure
of the tree we use the following conventions, in which
 different colors mean different ways of branching:
 \\

\tikzstyle{level 1}=[level distance=4cm, sibling distance=4cm]
\tikzstyle{level 2}=[level distance=5cm, sibling distance=3cm]

\begin{center}
\resizebox{0.55\textwidth}{!}{
\begin{tikzpicture}[scale=0.9,grow=right, sloped][H]
\node[rectangle,draw,purple,text=black,thick]{\shortstack{time\\position}}
    child {
        node[rectangle,draw,purple,text=black,thick]{\shortstack{time\\position}}
        child{
        node[rectangle,draw,thick]{...}
        edge from parent
        node[above]{label}
        node[below]{mark}
        }
        child{
        node[rectangle,draw,thick]{...}
        edge from parent
        node[above]{label}
        node[below]{mark}
        }
        edge from parent
        node[above]{label}
        node[below]{mark}
    }
    child {
        node[rectangle,draw,blue,text=black,thick]{\shortstack{time\\position}}
        edge from parent
        node[above]{label}
        node[below]{mark}
    };
\end{tikzpicture}
}
\end{center}

\noindent
Specifically, let us draw a tree sample for the following PDE:
$$ \frac{\partial u}{\partial t}
-\eta ( - \Delta / 2 ) u + c_0 (t,x) + c_{0,1}(t,x) u \frac{\partial u}{\partial x_1} = 0 $$
in dimension $d=1$. For this tree, there are two types of branching: we can either branch into no branch at all (which is represented in blue), or into two branches (one bearing the mark $0$ and one bearing the mark $1$, which is represented in purple).
The black color is used for leaves that have reached the horizon time $T$.
\\

\begin{center}
\resizebox{0.65\textwidth}{!}{
\begin{tikzpicture}[scale=0.9,grow=right, sloped]
\node[rectangle,draw,cyan,thick]{\shortstack{$t$\\$x$}}
    child {
        node[rectangle,draw,purple,text=black,thick] {\shortstack{$t+T_{\bar{1}}$\\ $X^{\bar{1}}_{T_{\bar{1}},x} $}}
            child {
                node[rectangle,draw,purple,text=black,thick] {\shortstack{$t+T_{\bar{1}} + T_{(1,2)}$\\ $X^{(1,2)}_{T_{(1,2)},x} $}}
                child{
                node[rectangle,draw,blue,text=black,thick]{\shortstack{$t+T_{\bar{1}}+ T_{(1,2)}+ T_{(1,2,2)}$\\ $X^{(1,2,2)}_{T_{(1,2,2)},x} $}}
                edge from parent
                node[above]{~~~$(1,2,2)$~~~}
                node[below]{$1$}
                }
                child{
                node[rectangle,draw,thick]{\shortstack{$T$\\ $X^{(1,2,1)}_{T,x} $}}
                edge from parent
                node[above]{$(1,2,1)$}
                node[below]{$0$}
                }
                edge from parent
                node[above] {$(1,2)$}
                node[below]  {$1$}
            }
            child {
                node[rectangle,draw,blue,text=black,thick] {\shortstack{$t+T_{\bar{1}} + T_{(1,1)}$\\ $X^{(1,1)}_{T_{(1,1)},x}  $}}
                edge from parent
                node[above] {$(1,1)$}
                node[below]  {$0$}
            }
            edge from parent
            node[above] {$\bar{1}$}
            node[below]{$0$}
    };
\end{tikzpicture}
}
\end{center}
 \noindent
In the above example we have
$\mathcal{K}^{\circ}= \{\bar{1}, (1,1) , (1,2) ,(1,2,2)\}$
and
 $\mathcal{K}^{\partial} = \{(1,2,1)\}$.

\section{Probabilistic representation}
\label{s3}
Given $t \in [0,T]$, $x\in \mathbb{R}^d$ and a mark $i\in \{0,1,\ldots ,d \}$,
we consider the functional $\mathcal{H}_\phi$
  of the random tree $\mathcal{T}_{t,x,i}$, defined as
\begin{equation}
\nonumber 
        {\mathcal{H}_\phi (\mathcal{T}_{t,x,i}) :=
          \prod_{\widebar{k} \in \mathcal{K}^{\circ}} \frac{c_{I_{\widebar{k}}}\big(T_{\widebar{k}},X^{\widebar{k}}_{T_{\widebar{k}},x}\big)\mathcal{W}_{\widebar{k}}}{q_{I_{\widebar{k}}}\rho(\tau_{\widebar{k}})} \prod_{\widebar{k} \in \mathcal{K}^{\partial}} \frac{\big( \phi\big(X^{\widebar{k}}_{T,x}\big)-\phi\big(X^{\widebar{k}}_{T_{\widebar{k}-},x}\big){\mathbbm{1}_{\{\theta_{\widebar{k}} \neq 0\}}}\big) \mathcal{W}_{\widebar{k}}}{\widebar{F}(T-T_{\widebar{k}-})}},
\end{equation}
 where $\widebar{F}(z) : = 1 - \mathbb{P}(T_{\widebar{1}} \leq z )$, $z\geq 0$,
 and $\mathcal{W}_{\widebar{k}}$ is a random weight defined by

\begin{subequations}
\begin{empheq}[left=\hskip-0.cm\text{$\mathcal{W}_{\widebar{k}} :=$}\empheqlbrace]{align}
 \label{fjklsf0}
  &
  \displaystyle
\mathbbm{1}_{\{\theta_{\widebar{k}}=0\}}
+ \mathbbm{1}_{\{\theta_{\widebar{k}}\neq 0\}}
\frac{\big(X^{\widebar{k}}_{T_{\widebar{k}},x}
  -X^{\widebar{k}}_{T_{\widebar{k}-},x}\big)_{\theta_{\widebar{k}}}}{S^{\widebar{k}}_{T_{\widebar{k}}}-S^{\widebar{k}}_{T_{\widebar{k}-}}}
 & \mbox{if~}~ \widebar{k}\in \mathcal{K}^{\circ},
\\
\label{fjklsf} 
 &  \displaystyle
 \mathbbm{1}_{\{\theta_{\widebar{k}}=0\}} + \mathbbm{1}_{\{\theta_{\widebar{k}} \neq 0\}} \frac{\big(X^{\widebar{k}}_{T,x}-X^{\widebar{k}}_{T_{\widebar{k}-},x}\big)_{\theta_{\widebar{k}}}}{S^{\widebar{k}}_T-S^{\widebar{k}}_{T_{\widebar{k}-}}}
 & \mbox{if~}~ \widebar{k}\in \mathcal{K}^{\partial},
\end{empheq}
\end{subequations}
 where $\theta_{\widebar{k}} \in \{0,1,\ldots ,d \}$
denotes the mark of the particle $\widebar{k}$,
and
for each particle labeled
$\widebar{k} = (1,k_2,\ldots ,k_n) \in \mathbb{N}^n$
at generation $n\geq 1$,
the subordinator
$\big(S_t^{\widebar{k}}\big)_{t\in [T_{\widebar{k}-},T_{\widebar{k}}]}$
is defined as
$$
S_t^{\widebar{k}} := S_{t-T_{\widebar{k}-}}^{n,\pi_n(\widebar{k})},
\qquad t\in [T_{\widebar{k}-},T_{\widebar{k}}].
$$
Here,
$\big(X^{\widebar{k}}_{T,x}-X^{\widebar{k}}_{T_{\widebar{k}-},x}\big)_{\theta_{\widebar{k}}}$
denotes the $\theta_{\widebar{k}}$-$th$ component of the vector
$X^{\widebar{k}}_{T,x}-X^{\widebar{k}}_{T_{\widebar{k}-},x}$ in $\real^d$.

\begin{theorem}{}
\label{t1}
 Let $m\geq 0$ denote the number of partial derivatives
$\partial u(t,x)/\partial x_1, \ldots,\partial u(t,x)/\partial x_m$
appearing in \eqref{eq:1}, and let $m_0 \in \{m,\ldots , d\}$.
Under the integrability conditions
 \begin{equation}
   \label{eq:IC}
         {\mathbb{E}\left[
             \sup_{x\in \real^d}
             \big|
             \mathcal{H}_\phi (\mathcal{T}_{t,x,i})\big|
             \right] < \infty},
    \quad
    t\in [0,T],
    \quad
    i=0,1,\ldots ,m,
\end{equation}
 and
\begin{equation}
 \label{eq:IC2}
          {\mathbb{E}\left[
             \big|
             \mathcal{H}_\phi (\mathcal{T}_{t,x,i})\big|
             \right] < \infty},
          \quad
(t,x) \in [0,T]\times \real^d,
\quad
  i=m+1,\ldots ,m_0,
\end{equation}
the function
 \begin{equation}
   \label{u0tx2}
u (t,x) :=  \mathbb{E}\big[ \mathcal{H}_\phi (\mathcal{T}_{t,x,0})\big],
\quad (t,x) \in [0,T]\times \real^d,
\end{equation}
 is an {integral} solution of the PDE \eqref{eq:1}.
 In addition, the 
 partial derivatives $\partial u(t,x) / \partial x_i$ exist and are represented as
 \begin{equation}
   \label{u0tx2-1}
\frac{\partial u}{\partial x_i} (t,x) :=  \mathbb{E}\big[ \mathcal{H}_\phi (\mathcal{T}_{t,x,i})\big],
\quad (t,x) \in [0,T]\times \real^d,
\quad i=1,\ldots ,m_0.
\end{equation}
\end{theorem}
\begin{Proof}
  We denote by $\varphi (t,y-x)$ the kernel of the
  pseudo differential operator $-\eta(-\Delta / 2)$,
  which is the fundamental solution of the PDE
  $\partial \varphi / \partial t = -\eta(-\Delta / 2) \varphi$.
 Letting
$$
 u_i(t,x) :=
\mathbb{E}\big[ \mathcal{H}_\phi (\mathcal{T}_{t,x,i})\big],
\quad (t,x) \in [0,T]\times \real^d, \quad i =1,\ldots ,m_0,
$$
 and applying the Markov property at
  the first branching time $T_{\widebar{1}}$ on the tree
  $\mathcal{T}_{t,x,0}$, we have
    \begin{eqnarray}
    \nonumber
  u(t,x) & := & \mathbb{E}\big[ \mathcal{H}_\phi \big(\mathcal{T}_{t,x,0}\big)\big]
  \\
    \nonumber
  &= &  \mathbb{E}\big[ \mathcal{H}_\phi \big(\mathcal{T}_{t,x,0}\big)\big(
    \mathbbm{1}_{\{T_{\widebar{1}}>T-t \}} +\mathbbm{1}_{\{T_{\widebar{1}}\leq T-t \}}\big) \big]
  \\
    \nonumber
  &= & \frac{\mathbb{P}(T_{\widebar{1}}>T-t)}{\widebar{F}(T-t)} \int_{\real^d} \varphi (T-t,y-x) \phi(y) dy
  \\
    \nonumber
  & &
  + \sum_{l = (l_0,\ldots ,l_m) \in {\cal L}_m}
  \int_0^{T-t} \int_{\real^d} \varphi (s,y-x)
  c_l(t+s,y) u^{l_0}(t+s,y) \prod_{j=1}^m u_j^{l_j}(t+s,y)dy ds
  \\
    \nonumber
  &= & \int_{\real^d} \varphi (T-t,y-x) \phi(y) dy
  \\
  \label{fdsfas}
  & &
  +
  \sum_{ l = (l_0,\ldots ,l_m) \in {\cal L}_m}
  \int_t^T \int_{\real^d} \varphi (s-t,y-x)
  c_l(s,y) u^{l_0}(s,y)\prod_{j=1}^m u_j^{l_j} (s,y) dy ds.
\end{eqnarray}
  Next, if $\theta_{\widebar{1}} \in \{1,\ldots ,d\}$, we have
\begin{enumerate}[a)]
    \item the subordination relation
$$
 X_{s,x}^{\widebar{1}} = x + Z_{s-t}^{1,1}
 = x + B_{S_s} - B_{S_t} \stackrel{d}{=}
 x + B_{S_s- S_t}, \quad
 x\in \real^d, \quad 0 \leq t \leq s \leq T,
 $$
 where $(B_t)_{t\in \real_+}
 = ((B_t)_1,\ldots , (B_t)_d)_{t\in \real_+}$ is a standard $d$-dimensional
 Brownian motion,
\item  conditional integration by parts
 with respect the Gaussian density of the
 $\theta_{\widebar{1}}$-th component
 $(B_{S_s- S_t})_{\theta_{\widebar{1}}}$ given $S_s-S_t$, 
 and
 \item the definition \eqref{fjklsf0}-\eqref{fjklsf} of
   $\mathcal{W}_{\widebar{1}}$,
  \end{enumerate}
  we have
\begin{eqnarray}
\nonumber
  \mathbb{E}\big[
    \mathcal{W}_{\widebar{1}} h \big(X^{\widebar{1}}_{t,x} \big)
    \ \!  \big| \ \!  T_{\widebar{1}} > T-t \big]
& = &
\mathbb{E}\left[
  \frac{(B_{S_T - S_t})_{\theta_{\widebar{1}}}}{S_T-S_t} h \big(X^{\widebar{1}}_{t,x} \big)
\ \!  \Big| \ \!  T_{\widebar{1}} > T-t \right]
\\
\nonumber
& = &
\mathbb{E}\left[
   \frac{\partial h}{\partial x_{\theta_{\widebar{1}}}}
   \big( x + B_{S_T- S_t} \big)
 \right]
\\
\nonumber
&= &
\frac{\partial }{\partial x_{\theta_{\widebar{1}}}}
\int_{\real^d} \varphi (T-t,y)
h (x+y)dy
\\
  \label{fskl2}
 & = &
\frac{\partial }{\partial x_{\theta_{\widebar{1}}}}
\mathbb{E}\big[
h \big( x + B_{S_T- S_t} \big)
 \big],
\end{eqnarray}
 for any function $h$ in the space
 ${\cal C}^1_b (\real^d)$ of ${\cal C}^1$ bounded functions on $\real^d$.
 As in Theorem~3.1 in \cite{fournie},
 see the proof argument of Corollary~3.6
 in \cite{kawai}, the above identity \eqref{fskl2} extends from
 $h \in {\cal C}^1_b (\real^d)$
 to $\phi ( x + \cdot )$, 
 with $\phi : \real^d \to \real$ continuous and bounded,
 as the differentiability relation
\begin{equation}
  \label{fjkdsl3}
  \frac{\partial }{\partial x_{\theta_{\widebar{1}}}}
\int_{\real^d} \varphi (T-t,y)
 \phi (x+y) dy
 =
\mathbb{E}\big[
  \mathcal{W}_{\widebar{1}}
  \big( \phi \big(X^{\widebar{1}}_{t,x} \big) - \phi ( x ) \big)
   \ \!  \big| \ \!  T_{\widebar{1}} > T-t \big],
\end{equation}
 which holds from \eqref{fskl2}
 and the fact that $\mathbb{E}\big[\mathcal{W}_{\widebar{1}}]=0$.
 Next, noting that by \eqref{eq:IC} and dominated convergence,
the function
$$
h (y) := c_l(s,y) u^{l_0}(s,y)\prod_{j=1}^m u_j^{l_j} (s,y),
\qquad y\in \real^d,
$$
is continuous and bounded, a similar argument shows that
\begin{eqnarray}
  \nonumber
  \lefteqn{
  \frac{\partial }{\partial x_{\theta_{\widebar{1}}}}
  \int_{\real^d} \varphi (s-t,y)
  c_l(s,x+y) u^{l_0}(s,x+y) \prod_{j=1}^m u_j^{l_j}(s,x+y)
  dy
  }
  \\
  \nonumber
  & = &
  \frac{\partial }{\partial x_{\theta_{\widebar{1}}}}
  \mathbb{E}\left[
  c_l \big(s, x + B_{S_s- S_t} \big) u^{l_0} \big(s,x + B_{S_s- S_t} \big)
  \prod_{j=1}^m u_j^{l_j} \big(s,x + B_{S_s- S_t} \big)
 \right]
  \\
  \label{fskl2-2}
  & = &
  \mathbb{E}\left[
      \mathcal{W}_{\widebar{1}}
       c_l(s,X^{\widebar{1}}_{s,x}) u^{l_0}(s,X^{\widebar{1}}_{s,x})\prod_{j=1}^m u_j^{l_j} (s,X^{\widebar{1}}_{s,x})
    \ \!  \bigg| \ \!  T_{\widebar{1}}=s \right]
 ,
\end{eqnarray}
$l = (l_0,\ldots ,l_m) \in {\cal L}_m$,
$0\leq t \leq s \leq T$.
 Applying the Markov property at
  the first branching time $T_{\widebar{1}}$ on the tree
  $\mathcal{T}_{t,x,i}$ and using \eqref{eq:IC}-\eqref{eq:IC2}
  and \eqref{fjkdsl3} we have, for $i=1,\ldots ,m_0$,
\begin{eqnarray*}
u_i(t,x) &= & \mathbb{E}\big[ \mathcal{H}_\phi (\mathcal{T}_{t,x,i})
  \big(
  \mathbbm{1}_{\{T_{\widebar{1}}>T-t \}} +\mathbbm{1}_{\{T_{\widebar{1}}\leq T-t \}}\big) \big]
\\
&= & \frac{\mathbb{P}(T_{\widebar{1}} >T-t)}{\widebar{F}(T-t)} \mathbb{E}\big[ \mathcal{W}_{\widebar{1}} \big( \phi \big(X^{\widebar{1}}_{t,x} \big)
  - \phi ( x ) \big) \ \! \big| \ \! T_{\widebar{1}} > T-t \big]
\\
& & +
\sum_{l = (l_0,\ldots ,l_m) \in {\cal L}_m}
\int_t^T \mathbb{E}\Bigg[\frac{\mathcal{W}_{\widebar{1}}}{\rho (T_{\widebar{1}})} c_l(s,X^{\widebar{1}}_{s,x}) u^{l_0}\big(s,X^{\widebar{1}}_{s,x}\big)\prod_{j=1}^m u_j^{l_j} \big(s,X^{\widebar{1}}_{s,x}\big)
  \ \!  \bigg| \ \!  T_{\widebar{1}}=s \Bigg] \rho(s) ds\\
&= &
\frac{\partial }{\partial x_i}
\int_{\real^d} \varphi (T-t,y) \phi (x+y) dy 
\\
& &
+ \sum_{ l = (l_0,\ldots ,l_m) \in {\cal L}_m}
\frac{\partial }{\partial x_i}
\int_t^T \int_{\real^d} \varphi (s-t,y)
\left(
c_l(s,x+y) u^{l_0}(s,x+y) \prod_{j=1}^m u_j^{l_j}(s,x+y)
\right)
dy ds.
\end{eqnarray*}
By \eqref{fdsfas} this shows \eqref{u0tx2-1}, i.e.
$$ u_i (t,x) = \frac{\partial u}{\partial x_i}(t,x),
\quad (t,x)\in [0,T] \times \real^d, \quad i=1,\ldots ,m_0,
$$
 and therefore we have
 \begin{eqnarray*}
   u(t,x) & = & \int_{\real^d} \varphi (T-t,y-x)  \phi(y)  dy
   \\
   & &
   + \sum_{l = (l_0,\ldots ,l_m) \in {\cal L}_m}
   \int_t^T \int_{\real^d} \varphi (s-t,y-x)
   c_l(s,y) u^{l_0}(s,y) \prod_{j=1}^m \left(
   \frac{\partial u}{\partial y_j} (s,y)\right)^{l_j} dy ds,
  \end{eqnarray*}
 $(t,x)\in [0,T] \times \real^d$, 
 showing that $u$ is an {integral} solution of \eqref{eq:1}.
\end{Proof}
We note that \eqref{eq:IC} also implies that
$u(t,\cdot )$ and $u_i(t,\cdot )$ are in $L^\infty (\real^d)$ for all $t\in [0,T]$
and $i=1,\ldots ,m$.
{
In the next proposition, we note that the probabilistic representation of
Theorem~\ref{t1} can be used to recover 
the classical result of \cite{fujita} on the blow-up of semilinear PDEs,
in the case of the fractional Laplacian. 
\begin{prop}
  \label{fkdslf1}
   (\cite{fujita}, \cite{sugitani}, \cite{blmw})  
 Consider the PDE
  \begin{equation}
    \label{dvkcx}
    \frac{\partial u}{\partial t} + \Delta_{\alpha} u + u^{1+\beta}=0
  \end{equation}
  with strictly positive terminal condition $u(T,x)=\phi(x)>0$,
  $x\in \real^d$.
  Under Assumption~(\hyperlink{BGJhyp}{$A$}), 
  when $\alpha \geq \beta d$ 
  there exists $T>0$ such that \eqref{dvkcx} 
  admits no solution on $[0,T]$. 
\end{prop} 
\begin{Proof}
Given $\varphi $ the solution of the heat equation 
$\partial_t \varphi  + \Delta_{\alpha} \varphi =0$
with $\varphi (T,x)=\phi(x)$,
 we denote as $v(t,x;T)$ the unique solution of 
 \begin{equation}
 \nonumber 
   \partial_t v(t,x;T) + \Delta_{\alpha} v(t,x;T) + v(t,x;T) \varphi^\beta (t,x)=0,
   \qquad v(T,x;T)=\phi(x),
   \end{equation}
 $(t,x)\in [0,T]\times \real^d$,
 which is a sub-solution of \eqref{dvkcx}.
 Since $\varphi$ and $\phi$ are bounded on $\real^d$,
 $v(t,x;T)$ can be represented
 by Theorem~\ref{t1} using a $1$-branching tree as 
$$
v(t,x;T) =
\mathbb{E}\big[ \mathcal{H}_\phi (\mathcal{T}_{t,x,0})\big]
= 
\mathbb{E}_{t,x}\left[\prod_{\bar{k} \in \mathcal{K}^{\circ}}
  \frac{\varphi^\beta( T_{\bar{k}},X_{T_{\bar{k}}})}{\rho(\Delta T_{\bar{k}})}
  \prod_{\bar{k} \in \mathcal{K}^{\partial}}
  \frac{\phi(X_{T_{\bar{k}}})}{\widebar{F}(T-T_{\bar{k}})} \right],
$$
 $(t,x) \in [0,T]\times \real^d$, where
$\mathbb{E}_{t,x}$ denotes the conditional expectation given that
 the tree is rooted at $(t,x)$. 
Next, letting ${\cal B}_r$ denoting the ball of radius $r>0$ centered at $0$
in $\real^d$,
 consider the event 
 $$ A := \big\{\omega \in \Omega \ : \ X_{T_{\bar{k}}}^{\bar{k}} \in {\cal B}_{(T-T_{\bar{k}})^{1/\alpha}}, 
 \bar{k} \in \mathcal{K}^{\circ},
 \
 and
 \
 X_{T}^{\bar{k}} \in {\cal B}_1, \bar{k} \in \mathcal{K}^{\partial}
 \big\}.
$$
 Let $x\in {\cal B}_1$ and denote by 
 $\mathcal{G} := \sigma \big( \tau^{i,j}, \ i,j \geq 1 \big)$,
 the sigma-algebra generated by the branching times. 
 By Lemma~2.2 in \cite{blmw} there exists $\kappa >0$ such that 
 $\mathbb{P}_x ( A \mid \mathcal{G} ) > \kappa >0$,
 a.e. on the event
$$
 B(t) := \left\{\omega \in \Omega, t\leq T_{\bar{k}-} \leq T/2, \bar{k} \in \mathcal{K}^{\partial} \right\} \cup \{ T_1 \geq T \},
$$ 
 where for $\omega \in B(t)$, the random tree $\mathcal{T}_{t,x}(\omega)$ has its last branching time before $T/2$.
 By (2.3) in \cite{blmw}, there exists $c>0$ such that 
 $$f(t) : = c (T-t)^{-d/\alpha} \int_{{\cal B}_{(T-t)^{1/\alpha}}} \phi(y)dy \leq \varphi (t,x),
 \qquad x \in {\cal B}_{(T-t)^{1/\alpha}}.
$$
 Hence,
 letting $\displaystyle C := \inf_{x\in {\cal B}_1} \phi (x) > 0$,
 we have 
 \begin{equation*}
\begin{split}
  v(0,x;T) &= \mathbb{E}_{0,x}\left[
    \prod_{k \in \mathcal{K}^{\circ}} \frac{\varphi^\beta( T_{\bar{k}},X_{T_{\bar{k}}})}{\rho(\Delta T_{\bar{k}})}\prod_{k \in \mathcal{K}^{\partial}} \frac{\phi(X_{T_{\bar{k}}})}{\widebar{F}(T-T_{\bar{k}})}
    \right]\\
  &\geq \mathbb{E}_{0,x}\left[
    \mathbbm{1}_A \mathbbm{1}_{B(0)} \prod_{k \in \mathcal{K}^{\circ}} \frac{f^\beta( T_{\bar{k}})}{\rho(\Delta T_{\bar{k}})}\prod_{k \in \mathcal{K}^{\partial}} \frac{C}{\widebar{F}(T-T_{\bar{k}})} \right]\\
  &\geq \mathbb{E}_{0,x}\left[
    \P ( A \mid \mathcal{G} )
    \mathbbm{1}_{B(0)} \prod_{k \in \mathcal{K}^{\circ}} \frac{f^\beta( T_{\bar{k}})}{\rho(\Delta T_{\bar{k}})}\prod_{k \in \mathcal{K}^{\partial}} \frac{C}{\widebar{F}(T-T_{\bar{k}})}      \right]
\\
&\geq \kappa  \mathbb{E}_{0,x}\left[ \mathbbm{1}_{B(0)} \prod_{k \in \mathcal{K}^{\circ}} \frac{f^\beta( T_{\bar{k}})}{\rho(\Delta T_{\bar{k}})}\prod_{k \in \mathcal{K}^{\partial}} \frac{C}{\widebar{F}(T-T_{\bar{k}})} \right]
\\
& : = \kappa g(0;T), 
\end{split}{}
\end{equation*}
 where the function
 $$
 g(t;T):=
 \mathbb{E}_{t,x}\left[ \mathbbm{1}_{B(t)} \prod_{k \in \mathcal{K}^{\circ}} \frac{f^\beta( T_{\bar{k}})}{\rho(\Delta T_{\bar{k}})}\prod_{k \in \mathcal{K}^{\partial}} \frac{C}{\widebar{F}(T-T_{\bar{k}})} \right], \qquad 0 \leq t\leq T, 
 $$
 is the solution of the ODE 
 \begin{equation}
   \label{fjdskl34}
   g(t;T) = C + \int_t^{T/2} f^\beta(s) g(s;T) ds, \qquad
  0 \leq t \leq \frac{T}{2}. 
 \end{equation}
 After solving \eqref{fjdskl34}, we obtain 
\begin{eqnarray*}
  g(t;T) &= & C \exp \bigg( \int_t^{T/2} f^\beta(s) ds \bigg) 
  \\
  &\geq & C
  \exp \bigg( \int_{T/ 2}^{T-t} \frac{\big(\int_{{\cal B}_{s^{1/\alpha}}} \phi(x)dx\big)^\beta}{s^{\beta d / \alpha}} ds
  \bigg)
  \\
  &\geq &
  C \exp \left(
  \bigg(\int_{\real^d} \phi(x)dx\bigg)^\beta
    \frac{
   (T-t)^{1-\beta d / \alpha }-(T/2)^{1- \beta d / \alpha } 
  }{2^{\beta}(1-\beta d/ \alpha)}
  \right), 
\end{eqnarray*}
hence $\lim_{T\to \infty} g(0;T)= \infty$, provided that $\alpha > \beta d$. 
Therefore, we have
$$
\lim_{T\to \infty} \inf_{x\in {\cal B}_1}|v(0,x;T)| = \infty,
$$
which is sufficient to conclude to blow-up as in \S~3 of \cite{blmw}.
In the critical case $d = \beta / \alpha$ we find 
$$ 
g(t;T) \geq C
\left( 2 \frac{T-t}{T} \right)^{ 
  \left(\int_{\real^d} \phi(x)dx / 2 \right)^\beta
},
\qquad 0 \leq t \leq T. 
  $$
 Letting now $w$ denote the solution of
 $\partial_t w + \Delta_{\alpha} w + w v^\beta=0$,
 with $w(T,x;T)=\phi(x)$, $x\in \real^d$, the above argument shows that
 $w(0,x;T) \geq \kappa h(0;T)$,
 where
$$ 
  h(t;T) = C \exp \left(   C^\beta \int_t^{T/2}
\left( \frac{T-s}{T/2} \right)^{ 
  \beta \left(\int_{\real^d} \phi(x)dx / 2 \right)^\beta
} ds \right), 
$$
and $\lim_{T\to \infty} h(0;T) = \infty$, therefore
$\lim_{T\to \infty} \inf_{x\in {\cal B}_1}|w(0,x;T)| = \infty$, 
which allows us to conclude to blow-up as above. 
Finally, the blow-up of $u(t,x;T)$ follows from the inequalities
 $u(t,x;T) \geq v(t,x;T) \geq w(t,x;T)$,
$(t,x) \in [0,T]\times \real^d$. 
\end{Proof}
} 
\section{$L^p$ Integrability} 
\label{s4}
In Theorem~\ref{t2} {and Proposition~\ref{t3}}
we derive sufficient conditions for the
integrability conditions \eqref{eq:IC}-\eqref{eq:IC2} to hold.
 As in Theorem~\ref{t1}, we let
 $m\geq 0$ denote the number of partial derivatives
 appearing in \eqref{fjkldsf}.
 The next result covers the case $\alpha = 2$ of the standard Laplacian
 by taking $\eta(\lambda) :=  2\lambda$ with 
 the deterministic subordinator $S_t = t$, $t\in \real_+$.
\begin{theorem}{}
  \label{t2}
  Under Assumption~(\hyperlink{BGJhyp}{$A$}),
 for any $p\geq 1$ and $m_0\in \{m,\ldots , d\}$ there exists
  a small enough $T = T(p,m_0)>0$ 
  such that
  \begin{equation}
    \label{fdskl3}
         {\mathbb{E}\left[
             \sup_{x\in \real^d}
             \big|
             \mathcal{{H}}_\phi (\mathcal{T}_{t,x,i})\big|^p
             \right] < \infty},
    \quad
    t\in [0,T],
    \qquad
    i=0,\ldots ,m_0,
\end{equation}
  provided that 
  \begin{equation}
    \label{fjkdsl3-1} 
{
   \int_0^T
\frac{1}{\rho^{p-1}(s)}
ds < \infty
\quad
\mbox{and}
\quad
\int_0^T \int_0^\infty
\frac{e^{-s\eta(\lambda)}}{\rho^{p-1}(s)} \lambda^{p/2-1}
d\lambda ds
< \infty}. 
\end{equation}
 When $p=1$, {both conditions in \eqref{fjkdsl3-1} are satisfied if }
 \begin{equation}
  \label{cd}
  \int_{\lambda_0}^\infty
  \frac{1}{\eta(\lambda) \sqrt{\lambda} }
  d\lambda < \infty
\end{equation}
 for some $\lambda_0>0$. 
\end{theorem}
\begin{Proof}
  Under \eqref{cc1} and \eqref{cc2}, the random variable
  $\mathcal{{H}}_\phi (\mathcal{T}_{t,x,i})$ is bounded as
 \begin{eqnarray}
   \nonumber
   \big|
   \mathcal{{H}}_\phi (\mathcal{T}_{t,x,i})
   \big|
   & \leq & \prod_{\widebar{k} \in \mathcal{K}^{\circ}} \frac{|c_{I_{\widebar{k}}}|_{\infty}|\mathcal{W}_{\widebar{k}}|}{q_{I_{\widebar{k}}}\rho(\tau_{\widebar{k}})} \prod_{\widebar{k} \in \mathcal{K}^{\partial}} \frac{\big|\big( \phi\big(X^{\widebar{k}}_{T,x}\big)-\phi\big(X^{\widebar{k}}_{T_{\widebar{k}-},x}\big)\mathbbm{1}_{\{\theta_{\widebar{k}} \neq 0\}}\big) \mathcal{W}_{\widebar{k}}\big|}{\widebar{F}(T-T_{\widebar{k}-})}
   \\
\nonumber
       & \leq &
   \prod_{\widebar{k} \in \mathcal{K}^{\circ}} \frac{|c_{I_{\widebar{k}}}|_{\infty}|\mathcal{W}_{\widebar{k}}|}{q_{I_{\widebar{k}}}\rho(\tau_{\widebar{k}})}
   \prod_{\widebar{k} \in \mathcal{K}^{\partial}}
   \frac{
          L \big\Vert X^{\widebar{k}}_{T,x}-X^{\widebar{k}}_{T_{\widebar{k}-},x}
          \big\Vert
          \big| \mathcal{W}_{\widebar{k}}\big|
          \mathbbm{1}_{\{\theta_{\widebar{k}} \neq 0\}}
     +
     |\phi|_\infty
     \mathbbm{1}_{\{\theta_{\widebar{k}} = 0\}}
          }{\widebar{F}(T-T_{\widebar{k}-})}
   \\
\nonumber
       & = &
   \prod_{\widebar{k} \in \mathcal{K}^{\circ}} \frac{|c_{I_{\widebar{k}}}|_{\infty}|\mathcal{W}_{\widebar{k}}|}{q_{I_{\widebar{k}}}\rho(\tau_{\widebar{k}})}
   \prod_{\widebar{k} \in \mathcal{K}^{\partial}}
   \frac{
     L \big\Vert
     Z_{T-T_{\widebar{k}-}}^{n,\pi_n(\widebar{k})}
          \big\Vert
          \big| \mathcal{W}_{\widebar{k}}\big|
          \mathbbm{1}_{\{\theta_{\widebar{k}} \neq 0\}}
     +
     |\phi|_\infty
     \mathbbm{1}_{\{\theta_{\widebar{k}} = 0\}}
          }{\widebar{F}(T-T_{\widebar{k}-})},
   \qquad
     x \in \mathbb{R}^d,
     \\
     \label{dsdfdnf$}
 \end{eqnarray}
 $i=0,\ldots ,m_0$.
 By the Cauchy-Schwartz inequality and \eqref{fdjskldf}, \eqref{fjklsf},
 when $\theta_{\widebar{k}} \in \{1,\ldots , d\}$ we have
 \begin{eqnarray}
   \nonumber 
    \mathbb{E}\big[ \big\Vert
      Z_{T-T_{\widebar{k}-}}^{n,\pi_n(\widebar{k})}
      \mathcal{W}_{\widebar{k}} \big\Vert^p \big]
    & = &
    \mathbb{E}\left[ \big\Vert
      Z_{T-T_{\widebar{k}-}}^{n,\pi_n(\widebar{k})}
      \big\Vert^p
      \frac{\big| \big(
        Z_{T-T_{\widebar{k}-}}^{n,\pi_n(\widebar{k})}
        \big)_{\theta_{\widebar{k}}}\big|^p}{(S^{\widebar{k}}_T-S^{\widebar{k}}_{T_{\widebar{k}-}})^p}
\right]
    \\
   \nonumber 
    & \leq &
    \sqrt{\mathbb{E}\Bigg[\frac{\big\Vert
          Z_{T-T_{\widebar{k}-}}^{n,\pi_n(\widebar{k})}
          \big\Vert^{2p}}{\big( S_{T_{\widebar{k}}}-S_{T_{\widebar{k}-}} \big)^p}
        \Bigg]\mathbb{E}\Bigg[
        \frac{\big(
Z_{T-T_{\widebar{k}-}}^{n,\pi_n(\widebar{k})}\big)_{\theta_{\widebar{k}}}^{2p}}{
        \big( S_{T_{\widebar{k}}}-S_{T_{\widebar{k}-}} \big)^p
        }
        \Bigg]}
    \\
    \label{fhdkjslf}
    & = & M_p \sqrt{d},
\end{eqnarray}
 where
$M_p: = \mathbb{E}[|X|^p] =
2^p \Gamma (p + 1/2 )/\sqrt{\pi}
$ for $X \sim \mathcal{N}(0,1)$.
Hence, by conditional independence
  given $\mathcal{G} := \sigma \big( \tau^{i,j},I^{i,j}, i,j \geq 1 \big)$ 
  of the terms in the product over
  $\widebar{k} \in \mathcal{K}^{\circ} \cup \mathcal{K}^{\partial}$
  in \eqref{dsdfdnf$},
  for all marks $i \in \{ 0,\ldots ,m_0 \}$ and all $t\in [0,T]$,
  denoting by $\mathbb{E}_{t,i}[\ \! \cdot \ \! ]$ the expected value given the
  initial mark $i$ at time $t$, we have
$$
\mathbb{E}\left[
             \sup_{x\in \real^d}
             \big|
             \mathcal{{H}}_\phi (\mathcal{T}_{t,x,i})\big|^p
             \right]   
  \leq \mathbb{E}_{t,i}\left[
    \prod_{\widebar{k} \in \mathcal{K}^{\circ}} \frac{|c_{I_{\widebar{k}}}|_{\infty}^p|\mathcal{W}_{\widebar{k}}|^p}{q_{\min}^{p-1} q_{I_{\widebar{k}}} \rho^p (\tau_{\widebar{k}})} \prod_{\widebar{k} \in \mathcal{K}^{\partial}} \frac{C_{\partial,p}}{\widebar{F}^p(T-T_{\widebar{k}-})} \right]
$$
  with $C_{\partial,p} := \max \big\{|\phi|_{\infty}^p,
  M_p L^p \sqrt{d} \big\}$
    and $q_{\min} := \min_{l = (l_0,\ldots ,l_m) \in {\cal L}_m} q_l >0$.
 To show \eqref{fdskl3} we will derive a system of Volterra
integral equations
  and give sufficient conditions for this system to have a local solution.
Proceeding by conditioning
  on the first branching time $T_{\widebar{1}}$
  as in the proof of Theorem~\ref{t1},
  we note that the functions $(v_0,v_1,\ldots v_{m_0})$
 defined as
 $$ v_i(t) :=
 \mathbb{E}_{t,i}\left[
    \prod_{\widebar{k} \in \mathcal{K}^{\circ}} \frac{|c_{I_{\widebar{k}}}|_{\infty}^p|\mathcal{W}_{\widebar{k}}|^p}{q_{\min}^{p-1} q_{I_{\widebar{k}}} \rho^p (\tau_{\widebar{k}})} \prod_{\widebar{k} \in \mathcal{K}^{\partial}} \frac{C_{\partial,p}}{\widebar{F}^p (T-T_{\widebar{k}-})} \right]
,
  \qquad i = 0,1,\ldots ,m_0,
  $$
  solve a system of
  Volterra integral equations of the form:
$$
v_0(t) = \frac{C_{\partial , p}}{\widebar{F}^{p-1} (T-t)} +
\sum_{l = (l_0,\ldots ,l_m) \in {\cal L}_m}
\frac{|c_l|_{\infty}^p}{q_{\min}^{p-1}}
\int_t^T \frac{v_0^{l_0}(s)}{\rho^{p-1}(s-t)}
   \prod_{j=1}^m v_j^{l_j}(s)
  ds,
$$
 and
\begin{eqnarray*}
\lefteqn{ v_i(t) =
  \frac{C_{\partial ,p}}{\widebar{F}^{p-1}(T-t)}
  }
 \\
  & & +
  \sum_{l = (l_0,\ldots ,l_m) \in {\cal L}_m} |c_l|_{\infty}^p
  \int_t^T
  \frac{1}{\rho^{p-1} (s-t)}
  \mathbb{E}\left[
    \int_{\real^d} G(S_s-S_t,x,y)\frac{|y-x|^p}{(S_s-S_t)^p}dy\right]
    v_0^{l_0}(s)\prod_{j=1}^m
    v_j^{l_j}(s) ds
    \\
  &=  &
   \frac{C_{\partial ,p}}{\widebar{F}^{p-1}(T-t)} +
   M_p
   \sum_{l = (l_0,\ldots ,l_m) \in {\cal L}_m}
\frac{|c_l|_{\infty}^p}{q_{\min}^{p-1}}
\int_t^T \frac{\mathbb{E}\big[S_{s-t}^{-p/2}\big]}{\rho^{p-1} (s-t)}
v_0^{l_0}(s)\prod_{j=1}^m v_j^{l_j}(s)
  ds,
\end{eqnarray*}
 for the marks $i=1,\ldots ,m_0$, where
$G(S_s-S_t,x,y)$ denotes the standard Gaussian kernel with
variance $S_s-S_t$, $0 \leq t <s$.
 We have
$$ v_0(t) \leq
 \frac{C_{\partial , p}}{\widebar{F}^{p-1} (T-t)} +
\sum_{l = (l_0,\ldots ,l_m) \in {\cal L}_m}
\frac{|c_l|_{\infty}^p}{q_{\min}^{p-1}}
\int_t^T \frac{v^{|l|}(s)}{\rho^{p-1}(s-t)}
ds,
$$
and
$$ v_i(t) \leq
   \frac{C_{\partial ,p}}{\widebar{F}^{p-1}(T-t)} +
   M_p
   \sum_{l = (l_0,\ldots ,l_m) \in {\cal L}_m}
\frac{|c_l|_{\infty}^p}{q_{\min}^{p-1}}
\int_t^T \frac{\mathbb{E}\big[S_{s-t}^{-p/2}\big]}{\rho^{p-1} (s-t)}
v^{|l|}(s)
  ds,
$$
 for the marks $i=1,\ldots ,m_0$.
Letting $v(t) := \max_{0\leq i \leq m} v_i(t)$,
$t\in [0,T]$, this leads to the Volterra integral inequality
\begin{eqnarray}
  \nonumber 
  v(t) & \leq &
  \frac{C_{\partial ,p}}{\widebar{F}^{p-1}(T-t)}
   + \max \Bigg\{
\sum_{l = (l_0,\ldots ,l_m) \in {\cal L}_m}
\frac{|c_l|_{\infty}^p}{q_{\min}^{p-1}}
\int_t^T \frac{v^{|l|}(s)}{\rho^{p-1}(s-t)}
ds,
\\
\nonumber
& &
   \qquad M_p
   \sum_{l = (l_0,\ldots ,l_m) \in {\cal L}_m}
\frac{|c_l|_{\infty}^p}{q_{\min}^{p-1}}
\int_t^T \frac{\mathbb{E}\big[S_{s-t}^{-p/2}\big]}{\rho^{p-1} (s-t)}
v^{|l|}(s)
  ds\Bigg\}
    \\
    \label{fdksf}
    &\leq &
    \frac{C_{\partial ,p}}{\widebar{F}^{p-1}(T-t)}
    +
    \sum_{l = (l_0,\ldots ,l_m) \in {\cal L}_m}
\frac{|c_l|_{\infty}^p}{q_{\min}^{p-1}}
\int_t^T
\frac{v^{|l|}(s)}{\rho^{p-1}(s-t)}
\big(
1 + M_p
   \mathbb{E}\big[S_{s-t}^{-p/2}\big]
   \big)
    ds.
\end{eqnarray}
\noindent
Using the comparison theorem for Volterra integral equations
(see page~121 of \cite{miller}), the integral inequality
\eqref{fdksf} admits a local in time solution
 $v(t) := \max_{0\leq i \leq m} v_i(t)$,
 provided that the corresponding Volterra
integral equation admits a local {maximal} solution $v(t)$
which is finite on an interval of the form $(T_*, T] \supset [0,T]$,
implying \eqref{fdskl3}.

  \medskip

In order to ensure the existence of this
local {in time maximal} solution, by {Theorem~5.1 page~116} Theorem~1 page~87 of \cite{miller}
it suffices to check that conditions (H3), (H4) and (H7)
pages~86-87 and 99 in \cite{miller} are satisfied, i.e.  
  \begin{equation}
  \tag{H3}
\nonumber 
  \sup_{0 \leq t \leq T}
  \int_t^T
\frac{1}{\rho^{p-1}(s-t)}
ds < \infty
\quad
\mbox{and}
\quad
\sup_{0 \leq t \leq T}
\int_t^T
\frac{   \mathbb{E}\big[ S_{s-t}^{-p/2}\big]}{\rho^{p-1}(s-t)}
ds
< \infty,
\end{equation}
 
\begin{equation}
  \tag{H4}
  \nonumber 
 \lim_{t \rightarrow t_0} \sup_{g \in C([0,T],\real ) \atop |g|_\infty \leq b}
      \int_0^T
    \left|
    \frac{  g^{|l|}(s) \mathbbm{1}_{\{t<s\}}}{\rho^{p-1}(s-t)}
  \big(
1 + M_p
   \mathbb{E}\big[S_{s-t}^{-p/2}\big]
   \big)
   -
   \frac{  g^{|l|}(s)\mathbbm{1}_{\{t_0<s\}}}{\rho^{p-1}(s-t_0)}
  \big(
1 + M_p
   \mathbb{E}\big[S_{s-t_0}^{-p/2}\big]
   \big)
   \right|
   ds
   = 0,
\end{equation}
    $l \in \mathcal{L}$, and
\begin{equation}
  \tag{H7}
  \label{c3}
\lim_{h \rightarrow 0 } \int_{t-h}^{t} \frac{1+M_p \mathbb{E}\big[S_{s-t+h}^{-p/2}\big] }{\rho^{p-1}(s-t+h)} ds = 0, 
\end{equation}
uniformly in $t \in [0,T]$. Regarding (H3), using \eqref{id} and \eqref{id2}
 we have
$$
  \mathbb{E}\big[S_{s-t}^{-p/2}\big]
    =
  \frac{1}{\Gamma(p/2)}
  \int_0^\infty e^{-(s-t)\eta(\lambda)} \lambda^{p/2-1} d\lambda,
$$ 
which shows by \eqref{fjkdsl3-1} that (H3) is satisfied.
 Regarding (H4), under the condition $|g|_{\infty}\leq b$, $g \in C([0,T],\real )$, we have
\begin{align}
\nonumber 
& \lim_{t \rightarrow t_0} \sup_{g \in C([0,T],\real ) \atop |g|_\infty \leq b}
    \left|
  \int_0^T
    \left(
    \frac{  g^{|l|}(s) \mathbbm{1}_{\{t<s\}}}{\rho^{p-1}(s-t)}
  \big(
1 + M_p
   \mathbb{E}\big[S_{s-t}^{-p/2}\big]
   \big)
   -
    \frac{  g^{|l|}(s)\mathbbm{1}_{\{t_0<s\}}}{\rho^{p-1}(s-t_0)}
  \big(
1 + M_p
   \mathbb{E}\big[S_{s-t_0}^{-p/2}\big]
   \big)
   \right)
   ds
   \right|
   \\
\nonumber 
   & \leq b^{|l|} \lim_{t \rightarrow t_0}
     \int_0^T
   \left|
    \frac{\mathbbm{1}_{\{t<s\}}}{\rho^{p-1}(s-t)}
-    \frac{\mathbbm{1}_{\{t_0<s\}}}{\rho^{p-1}(s-t_0)}
\right|
ds
   \\
\nonumber 
   & \quad +
b^{|l|} M_p \lim_{t \rightarrow t_0}
    \int_0^T
  \left|\mathbbm{1}_{\{t<s\}}
    \frac{\mathbb{E}\big[S_{s-t}^{-p/2}\big]}{\rho^{p-1}(s-t)}
        -
    \mathbbm{1}_{\{t_0<s\}}
    \frac{   \mathbb{E}\big[S_{s-t_0}^{-p/2}\big]}{\rho^{p-1}(s-t_0)}
    \right|
   ds
   \\
   \nonumber
   & = 
   \frac{   b^{|l|} M_p }{\Gamma(p/2)}
   \lim_{t \rightarrow t_0}
  \int_0^T
  \int_0^\infty
  \left|\mathbbm{1}_{\{t<s\}}
  \frac{ e^{-(s-t)\eta(\lambda)}}{\rho^{p-1}(s-t)}
     -
 \mathbbm{1}_{\{t_0<s\}} \frac{ e^{-(s-t_0)\eta(\lambda)} }{\rho^{p-1}(s-t_0)}
    \right| \lambda^{p/2-1} d\lambda
    ds
    \\
    \nonumber
    & =  0
\end{align} 
for all $l = (l_0,\ldots ,l_m) \in {\cal L}_m$
by Scheff\'e's lemma since by \eqref{fjkdsl3-1} and dominated convergence we have 
$$
 \lim_{t \rightarrow t_0}
     \int_0^T
       \frac{\mathbbm{1}_{\{t<s\}}}{\rho^{p-1}(s-t)}
ds
=
     \int_0^T
    \frac{\mathbbm{1}_{\{t_0<s\}}}{\rho^{p-1}(s-t_0)}
ds
$$
and
$$
 \lim_{t \rightarrow t_0}
  \int_0^T
  \int_0^\infty
  \mathbbm{1}_{\{t<s\}}
  \frac{ e^{-(s-t)\eta(\lambda)}}{\rho^{p-1}(s-t)}
   \lambda^{p/2-1} d\lambda
    ds
    =
      \int_0^T
  \int_0^\infty
  \mathbbm{1}_{\{t_0<s\}}
  \frac{ e^{-(s-t_0)\eta(\lambda)}}{\rho^{p-1}(s-t_0)}
   \lambda^{p/2-1} d\lambda
    ds. 
$$    
{Regarding \eqref{c3}, by \eqref{fjkdsl3-1} we have}
\begin{equation}
  {\lim_{h\rightarrow 0}\int_{t-h}^{t} \frac{1+M_p \mathbb{E}\big[S_{s-t+h}^{-p/2}\big] }{\rho^{p-1}(s-t+h)} ds = \lim_{h\rightarrow 0}\int_0^h \frac{ds}{\rho^{p-1}(s)} +
    \lim_{h\rightarrow 0}
    \frac{M_p}{\Gamma(p/2)}\int_0^h \int_0^\infty \frac{e^{-s \eta(\lambda)}}{\rho^{p-1}(s)}\lambda^{p/2-1} d\lambda ds = 0.}
\end{equation}
 When $p=1$ we have
\begin{eqnarray*}
 \int_0^{\infty} \lambda^{-1/2} \int_t^T e^{-(s -t) \eta(\lambda)} ds d\lambda
 &= &
 \int_0^{\infty} \frac{1-e^{-(T-t) \eta(\lambda)}}{\eta(\lambda)} \lambda^{-1/2} d\lambda, 
\end{eqnarray*}
 and we conclude from the facts that
 the integrand
$(1-e^{-(T-t) \zeta \eta(\lambda)})\lambda^{-1/2}/\eta(\lambda)$
 is equivalent to $(T-t)/\sqrt{\lambda}$ 
 as  $\lambda \to 0$, and to $\lambda^{-1/2} / \eta(\lambda) $
as $\lambda \to +\infty$. 
\end{Proof}
 The probabilistic representation \eqref{u0tx2}
 provided in Theorem~\ref{t1} will be used to
estimate the solution of \eqref{eq:1}
by Monte Carlo simulations in Section~\ref{s5}.
Finiteness of the second moment
of the functional $\mathcal{H}_\phi (\mathcal{T}_{t,x,i})$
 is needed in order to control the convergence via
the central limit theorem, and is ensured by
the sufficient conditions on $\rho$ and $\eta$
in Theorem~\ref{t2}.
{
  \begin{remark} 
  When $m=0$, i.e. the PDE
 \eqref{eq:1} does not contain any partial derivative,   
 it follows, by inspection of its proof, that Theorem~\ref{t2}
 holds by replacing
  \eqref{fjkdsl3-1} with the single condition 
$
 \displaystyle   \int_0^T
\frac{1}{\rho^{p-1}(s)}
ds 
< \infty$. 
\end{remark} 
}
\noindent
   When $p=1,2$ and $\eta(\lambda) = (2\lambda)^{\alpha / 2}$
the integrability condition \eqref{cd}
can be made more specific in the case of fractional Laplacians.
\begin{corollary}
  Consider the case $\eta(\lambda) =  (2\lambda)^{\alpha / 2}$
  of the fractional Laplacian $\Delta_\alpha = -(-\Delta)^{\alpha/2}$.
  \begin{enumerate}[i)]
  \item When $p=1$, the integrability conditions \eqref{eq:IC}-\eqref{eq:IC2}
  hold whenever $\alpha \in (1,2)$.
\item
 When $p=2$ and $\rho: \mathbb{R}^+ \rightarrow (0,\infty )$
is the gamma probability density function
$\rho(s) := s^{\delta-1} e^{-s} / \Gamma(\delta)$
for $\delta >0$,
the integrability conditions \eqref{eq:IC}-\eqref{eq:IC2} hold
 provided that
  $\delta < 2- 2 / \alpha$.
\end{enumerate}
\end{corollary}
\begin{Proof}
  \noindent
  $(i)$ When $p=1$, by \eqref{cd} it suffices to note that
  the function $1 / ( \lambda^{\alpha/2} \sqrt{\lambda})$
  is integrable at $+\infty$ if and only if $1/2+\alpha /2 > 1$.

\noindent
$(ii)$ When $p=2$ we have
{
  $   \int_0^T
 ds / \rho (s) < \infty
$
since $\delta < 2$,
and 
}
\begin{eqnarray*}
  \int_0^T
  \frac{1}{\rho(s)}
  \int_0^\infty
  e^{-s (2\lambda)^{\alpha / 2}}
  d\lambda
  ds
 & = &
  \frac{2}{\alpha}
  \int_0^T
  \frac{s^{-2/\alpha} }{\rho(s)}
  \int_0^\infty
  e^{- \mu }
  \mu^{-1+2/\alpha }
  d\mu
  ds
  \\
   & = &
  \frac{2}{\alpha}
 \Gamma(\delta)
  \Gamma (2/\alpha ) \int_0^T
 s^{1-\delta -2/\alpha } e^{s}
 ds
  \\
   & < &
 \infty,
\end{eqnarray*}
which holds since $\delta < 2 - 2/\alpha$,
hence \eqref{fjkdsl3-1} is satisfied.
\end{Proof}
{
 In the case of the fractional Laplacian,
 quantitative bounds on the time $T$
 satisfying \eqref{fdskl3}  
 and ensuring existence of solutions on $[0,T]$ by Theorem~\ref{t1},
 are derived in the next result. 
 Note that \eqref{vlkcxnv}-\eqref{vncs} 
 hold respectively for the gamma probability density function
$\rho(s) := s^{\delta-1} e^{-s} / \Gamma(\delta)$
 when
 $0< \delta < 1 - 1/\alpha$, resp. 
 $0< \delta < 1 - p/(\alpha (p-1))$ if $p>1$. 
\begin{prop}
\label{t3}
Let $p\geq 1$.
 Under Assumption~(\hyperlink{BGJhyp}{$A$}),
assume that $\eta(\lambda) = (2\lambda)^{\alpha / 2}$ 
 with $\alpha \in (1,2]$, 
let $q_{\min} := \min_{l = (l_0,\ldots ,l_m) \in {\cal L}_m} q_l >0$,
 $C_{\partial,p} := \max \big\{|\phi|_{\infty},
  (2L)^p \Gamma (p + 1/2 ) \sqrt{d/\pi} \big\}$, 
 and $m_0 \in \{ m,\ldots , d\}$. 
 Then, the bound \eqref{fdskl3} holds for all $t\in [0,T]$,
 provided that $T$ satisfies condition $(a)$ or condition $(b)$ below. 
\begin{enumerate}[a)]
\item
  The time $T$ is small enough so that
\begin{equation}
   \label{vlkcxnv} 
   C_{\circ,p}(T) :=  \frac{1}{q_{\min}^p } \sup_{l \in \mathcal{L}} |c_l|^p
    \max\left\{
\frac{2 \Gamma(p/\alpha)}{2^{p/2}\alpha \Gamma(p/2)} \sup_{s\in (0,T]}
\frac{s^{-p/\alpha} }{\rho^p(s)} ,  \sup_{s \in (0,T]} \frac{1}{\rho^p(s)} \right\} 
  \leq 1, 
\end{equation}
 and 
 \begin{equation}
   \label{vncs} 
  \frac{C_{\partial,p}}{\widebar{F}^p(T)} \leq 1. 
 \end{equation}
\item
  The time $T$ is small enough so that
\begin{equation}
\nonumber 
\widetilde{C}_{\circ,p} (T):=  \frac{1}{q_{\min}^{p-1} }
\sup_{l \in \mathcal{L}} |c_l|^{p-1}
\max \left\{ \frac{2 \Gamma(p/\alpha)}{2^{p/2}\alpha \Gamma(p/2)}
\sup_{s\in (0,T]}  \frac{s^{-p/\alpha} }{\rho^{p-1}(s)} , 
      \sup_{s \in (0,T]} \frac{1}{\rho^{p-1}(s)} \right\} 
  <\infty, 
\end{equation}
and 
\begin{equation}
   \label{abc}
   T  < \frac{1}{\widetilde{C}_{\circ,p}(T)} \int_{C_{\partial,p} / \widebar{F}^{p-1}(T)}^\infty  \left(
    \sum_{l \in \mathcal{L}} |c_l|_{\infty}x^{|l|} \right)^{-1} dx. 
\end{equation} 
\end{enumerate}
\end{prop}
\begin{Proof}
  $a)$ 
  By \eqref{fhdkjslf}
  and conditional independence
  of the terms in the product over
  $\widebar{k} \in \mathcal{K}^{\circ} \cup \mathcal{K}^{\partial}$ 
  in \eqref{dsdfdnf$} 
  given $\mathcal{G} := \sigma \big( \tau^{i,j},I^{i,j}, i,j \geq 1 \big)$, 
  denoting by $\mathbb{E}_{t,i}[ \ \! \cdot \ \! ]$ the expected value given the
  starting time $t\in [0,T]$ of the tree with initial
  mark $i \in i \in \{ 0,\ldots ,m_0 \}$, we have 
\begin{eqnarray}
  \nonumber
  \lefteqn{ 
\mathbb{E}\left[ \sup_{x \in \mathbb{R}^d}
             \big|
             \mathcal{{H}}_\phi (\mathcal{T}_{t,x,i})\big|^p 
             \right]   
}
\\
\nonumber
& \leq & \mathbb{E}_{t,i} \left[ \mathbb{E} \left[  \prod_{\widebar{k} \in \mathcal{K}^{\circ}} \frac{|c_{I_{\widebar{k}}}|_{\infty}^p|\mathcal{W}_{\widebar{k}}|^p}{q_{I_{\widebar{k}}}^p\rho^p(T_{\widebar{k}}-T_{\widebar{k}-})} \ \! \bigg| \ \! \mathcal{G} \right] \mathbb{E} \left [
   \prod_{\widebar{k} \in \mathcal{K}^{\partial}}
   \frac{
     L^p \big\Vert
     Z_{T_{\widebar{k}}-T_{\widebar{k}-}}^{n,\pi_n(\widebar{k})}
          \big\Vert^p
          \big| \mathcal{W}_{\widebar{k}}\big|^p
          \mathbbm{1}_{\{\theta_{\widebar{k}} \neq 0\}}
     + |\phi|_\infty ^p
     \mathbbm{1}_{\{\theta_{\widebar{k}} = 0\}}
   }{\widebar{F}^p(T-T_{\widebar{k}-})} \ \! \bigg| \ \! \mathcal{G} \right] \right]
\\
\label{fhkjsdf}
& \leq & \mathbb{E}_{t,i}\left[ \mathbb{E}\left[
    \prod_{\widebar{k} \in \mathcal{K}^{\circ}} \frac{|c_{I_{\widebar{k}}}|_{\infty}^p|\mathcal{W}_{\widebar{k}}|^p}{q_{\min}^{p-1} q_{I_{\widebar{k}}} \rho^p (T_{\widebar{k}}-T_{\widebar{k}-})} \ \! \bigg| \ \! \mathcal{G} \right]\prod_{\widebar{k} \in \mathcal{K}^{\partial}} \frac{C_{\partial,p}}{\widebar{F}^p(T)} \right]. 
\end{eqnarray}
Next, for a particle labeled $\bar{k}$ with mark $\theta_{\bar{k}} \neq 0$,
using \eqref{id3} and
\eqref{vlkcxnv} we have
\begin{eqnarray*} 
  \frac{|c_{I_{\widebar{k}}}|_{\infty}^p}{q_{\min}^{p-1}q_{I_{\widebar{k}}}}
  \mathbb{E} \left[ \frac{|\mathcal{W}_{\widebar{k}}|^p}{\rho^p (T_{\widebar{k}}-T_{\widebar{k}-})} \ \! \bigg| \ \! \mathcal{G} \right]
  & = &
  \frac{|c_{I_{\widebar{k}}}|_{\infty}^p}{q_{\min}^{p-1}q_{I_{\widebar{k}}}}
  \mathbb{E} \left[ \mathbb{E} \left[ \frac{|\mathcal{W}_{\widebar{k}}|^p}{ \rho^p (T_{\widebar{k}}-T_{\widebar{k}-})} \ \! \bigg| \ \! \sigma ( S^{i,j})\right] \ \! \bigg| \ \! \mathcal{G} \right]\\
  \\
  & = &
  \frac{|c_{I_{\widebar{k}}}|_{\infty}^p}{q_{\min}^{p-1}q_{I_{\widebar{k}}}}
  \mathbb{E} \left[\frac{M_p}{S_{T_{\widebar{k}}-T_{\widebar{k}-}}^{p/2}} \frac{1}{ \rho^p (T_{\widebar{k}}-T_{\widebar{k}-})} \ \! \bigg| \ \! \mathcal{G}\right]
  \\
  & = &
  \frac{2^{1-p}|c_{I_{\widebar{k}}}|_{\infty}^p  M_p \Gamma(p/\alpha)}{\alpha q_{\min}^{p-1}q_{I_{\widebar{k}}}\Gamma(p/2)}
  \mathbb{E} \left[
       \frac{(T_{\widebar{k}}-T_{\widebar{k}-})^{-p/\alpha}}{ \rho^p (T_{\widebar{k}}-T_{\widebar{k}-})} \ \! \bigg| \ \! \mathcal{G}\right]
  \\
   & \leq & C_{\circ,p}(T)
\end{eqnarray*}
 and 
\begin{equation*}
\begin{split}
  \frac{|c_{I_{\widebar{k}}}|_{\infty}^p}{q_{\min}^{p-1} q_{I_{\widebar{k}}} }
  \mathbb{E} \left[ \frac{1}{\rho^p (T_{\widebar{k}}-T_{\widebar{k}-})} \ \! \bigg| \ \! \mathcal{G}\right] \leq C_{\circ,p}(T), 
\end{split}
\end{equation*}
 hence by \eqref{fjklsf0}, under 
\eqref{vncs} the random variable $\mathcal{H}(\mathcal{T}_{t,x,i})$ is bounded by $1$.
\\
$b)$
 We rewrite \eqref{fhkjsdf} as 
\begin{equation*}
\begin{split}
\mathbb{E}\left[ \sup_{x \in \mathbb{R}^d}
             \big|
             \mathcal{{H}}_\phi (\mathcal{T}_{t,x,i})\big|^p 
             \right]   
   &\leq \eta(t) := \mathbb{E}_{t,i} \left[ \prod_{\widebar{k} \in \mathcal{K}^{\circ}} \frac{\widetilde{C}_{\circ,p}(T) | c_{I_{\bar{k}} } |_\infty}{q_{I_{\widebar{k}}}\rho(T_{\widebar{k}}-T_{\widebar{k}-})} 
   \prod_{\widebar{k} \in \mathcal{K}^{\partial}}
   \frac{ C_{\partial,p} / \widebar{F}^{p-1}(T)}{\widebar{F}(T-T_{\widebar{k}-})}\right] 
\end{split}
\end{equation*}
 where $\eta (t)$ solves the ODE 
$$
 \eta(t) = \frac{C_{\partial,p}}{ \widebar{F}^{p-1}(T)} + \widetilde{C}_{\circ,p}(T)
 \int_t^T
 \sum_{l \in \mathcal{L}}
 | c_l |_\infty
 \eta(s)^{|l|} ds, \qquad t\in [0,T], 
$$
which admits a (finite) solution as long as \eqref{abc} holds. 
\end{Proof} 
}
\subsubsection*{Examples}

We discuss some examples of subordinators and
their Laplace exponents in relation with
the above integrability conditions, see e.g. \S~6 of \cite{kyprianou}.

\medskip

\noindent
The first example is a variation of the stable subordinator.
\begin{example}{Sum of independant stable processes.}
For $a,b > 0$ and $0<\alpha<\beta<1$, let
$$ \eta(\lambda) := a \lambda^{\beta-\alpha} + b\lambda^{\beta},
$$
 which is the Laplace exponent of the sum of two independant stable subordinators with parameters $\beta - \alpha$ and $\beta$.
 Since we have $\eta(\lambda)\sqrt{\lambda} \sim b \lambda^{\beta+1/2}$
 as $\lambda$ tends to infinity,
 the integrability condition \eqref{cd} holds if and only if
$\beta > 1/2$.
\end{example}{}

\begin{example}{Stable subordinator with drift.}
  The Bernstein function
$$ \eta(\lambda) := \kappa + \mu \lambda + c\lambda^{\alpha}
$$
  is the Laplace exponent of an $\alpha$-stable subordinator,
  $\alpha \in (0,1)$, with drift $\mu > 0$ killed at the rate
  $\kappa>0$, with $c>0$.
  Due to the equivalent
  $\eta(\lambda)\sqrt{\lambda} \sim \lambda^{3 / 2}$
  as $\lambda$ tends to infinity,
  the integrability condition \eqref{cd} is always satisfied in this case.
\end{example}{}

\begin{example}
Consider the Bernstein function
$$ \eta(\lambda) := \frac{c\lambda \Gamma(\lambda + \nu)}{\Gamma(\lambda+\nu+\mu)} = \frac{c\lambda}{\Gamma(\mu)} B(\lambda+\nu,\mu),
$$
with $c>0$, $\nu \geq 0$, $\mu \in (0,1)$.
Due to the equivalent
$\sqrt{\lambda}\eta(\lambda) \sim c\lambda^{-\mu+3 / 2}$
as $\lambda$ tends to infinity,
the integrability condition \eqref{cd} holds if and only if
$\mu < 1/2$.
\end{example}{}

\begin{example}{Relativistic stable subordinator.}
  The Bernstein function
  $\eta(\lambda) := (\lambda+m^{2/\alpha})^{\alpha/2} - m$,
  with $\alpha \in (0,2)$, $m>0$,
  satisfies
  $\eta(\lambda)\sqrt{\lambda} \sim \lambda^{(1+\alpha )/2}$
  as $\lambda$ tends to infinity,
  thus the integrability condition \eqref{cd} holds if and only if
  $ \alpha > 1$.
\end{example}{}
\begin{example}
  For $\alpha \in (0,2)$ and $\beta \in (0,2-\alpha)$,
  the Bernstein function $\eta(\lambda) := \lambda^{\alpha/2}(\log(1+\lambda))^{\beta/2}$ satisfies the integrability condition
  \eqref{cd} if and only if $\alpha>1$.
 When $\beta \in (0,\alpha)$, the
Bernstein function $\eta(\lambda) := \lambda^{\alpha/2}(\log(1+\lambda))^{-\beta/2}$ satisfies the integrability condition \eqref{cd} if and only if $\alpha>1$.
\end{example}
The following table summarizes the above examples of integrability conditions.
\bigskip
\begin{center}
 \begin{tabular}{||c|c|c||}
 \hline
 Laplace exponent $\eta(\lambda)$ & Parameters & Integrability condition \\ [0.5ex]
 \hline\hline
$ a \lambda^{\beta-\alpha} + b\lambda^{\beta}$ &  $a,b > 0$ and $0<\alpha<\beta<1$ & $0 < \max (\alpha , 1/2 ) < \beta <1$\\
 \hline
 $ \kappa + \mu \lambda + c\lambda^{\alpha}$ & $\alpha \in (0,1)$, $\mu > 0$,
 $\kappa , c>0$ & Always satisfied \\
 \hline
$c\lambda B(\lambda+\nu,\mu) /\Gamma(\mu)$ & $c>0,~ \nu \geq 0,~ \mu \in (0,1)$ & $ 0 < \mu < 1 / 2 $\\  \hline
$(\lambda+m^{2/\alpha})^{\alpha/2} - m$ & $\alpha \in (0,2), ~m>0$ & $1 < \alpha < 2$\\
 \hline
 $\lambda^{\alpha/2}(\log(1+\lambda))^{\beta/2}$ & $\alpha \in (0,2), ~\beta \in (0,2-\alpha)$ & $1 < \alpha < 2$\\
 \hline
 $\lambda^{\alpha/2}(\log(1+\lambda))^{-\beta/2}$ & $\alpha \in (0,2), ~\beta \in (0,\alpha)$ & $1 < \alpha < 2$ \\
\hline
\end{tabular}
\end{center}
\subsubsection*{Higher order derivatives}
Here, we shortly discuss the difficulties in
dealing with higher orders of derivation
inside the coefficient $f$ of \eqref{eq:1}.
Writing the iterated integrations by parts relation
\eqref{fskl2}-\eqref{fskl2-2} for a higher order of derivation $p\geq 2$
would require to use a weight $\mathcal{W}_k$ given from a
Hermite polynomial of degree $p$, and therefore
to show the integrability of
$\big(B_{S_s - S_t} \big)^p / (S_s-S_t)^p$.
Since
$B_{S_s - S_t} / (S_s-S_t)^{1/2} \sim {\cal N}(0,1)$
given $S_s - S_t$, this would however require to show
the finiteness of
$$
 \int_0^T \mathbb{E}\big[ S_s^{-p/2} \big] ds
$$
 for $p \geq 2$, which does not hold.
 Indeed, from \eqref{id2}, we have
$$
  \int_0^T  \mathbb{E}\big[ S_s^{-p/2} \big] ds
  = \frac{1}{\Gamma(p/2)}
  \int_0^T \int_0^\infty e^{-s \eta(\lambda)} \lambda^{p/2-1} d\lambda ds
  = \frac{1}{\Gamma(p/2)}\int_0^\infty \frac{1-e^{-T \eta(\lambda)}}{\eta(\lambda)} \lambda^{p/2-1} d\lambda
  ,
$$
which is not integrable at $0$ when $p\geq 2$.
For example, in the case of the fractional Laplacian
when $(S_t)_{t\in \real_+}$ is an $\alpha / 2$-subordinator,
\eqref{id3} shows that
$$
\mathbb{E}\big[ S_s^{-p/2} \big]
= \frac{2^{1-p/2} \Gamma(p/\alpha)}{\alpha s^{p/\alpha} \Gamma(p/2)}
$$
 is integrable in $s$ around $0$ if and only if
$\alpha \in (p,2)$,
which excludes integration by parts of order $p\geq 2$.
As a result,
this method does not allow for higher order integration by parts,
and therefore it does not extend to the treatment of
higher order derivatives in the PDE \eqref{eq:1}.

\section{Numerical examples}
\label{s5}

In this section we consider numerical examples involving the fractional Laplacian $\Delta_\alpha $ and
 the $\alpha / 2$-stable subordinator
$(S_t)_{t\in \real_+}$
 with Laplace exponent $\eta(\lambda) = (2\lambda )^{\alpha / 2}$
 for $\alpha \in (1,2)$.
 For the generation of random samples of $S_t$, we use the formula

$$
\widetilde{S}_t := 2t^{2/\alpha}\frac{\sin( \alpha \big(U+\pi / 2)/2\big)}{\cos^{2/\alpha} (U)} \left(
\frac{\cos\big(U-\alpha(U+\pi / 2)/2\big)}{E}\right)^{-1+2/\alpha }
$$
based on the
Chambers-Mallows-Stuck (CMS) method,
where $U \sim U (-\pi / 2,\pi /2)$,
and $E \sim {\rm Exp} (1)$,
see Relation~(3.2) in \cite{weron1996chambers},
where $\psi_S (\lambda)$
denotes the L\'evy symbol
of $(S_t)_{t\in \real_+}$, see \eqref{ls}.
We start by testing our algorithm on an equation admitting a known solution.
For $k\geq 0$, we consider the function
$$
\Phi_{k,\alpha} (x) := (1-\Vert x\Vert^2)^{k+\alpha / 2}_+,
\qquad x\in \real^d,
$$
 which is Lipschitz if $k>1-\alpha/2$, and solves the Poisson problem
 $\Delta_\alpha \Phi_{k,\alpha} = -\Psi_{k,\alpha}$
 on $\real^d$, with
 \begin{align}
\nonumber
&  \Psi_{k,\alpha} (x)
\\
\label{pp}
& := \left\{
\begin{array}{ll}
  \displaystyle
 \frac{\Gamma( ( d+\alpha ) / 2)
   \Gamma(k+1+ \alpha / 2 )}{
   2^{-\alpha}
   \Gamma(k+1)\Gamma( d / 2 )}
~{}_2F_1\left(
  \frac{d+\alpha}{2},-k;\frac{d}{2};\Vert x\Vert^2
  \right), ~~ \Vert x\Vert\leq 1
  \medskip
  \\
\displaystyle
 \frac{2^{\alpha}\Gamma( ( d+\alpha ) /2 )
  \Gamma(k+1+ \alpha / 2 )}{\Gamma(k+1+ ( d+\alpha ) / 2 )
   \Gamma(- \alpha / 2 )
   \Vert x\Vert^{d+\alpha}
 }
  {}_2F_1\left(
  \frac{d+\alpha}{2},\frac{2+\alpha}{2};k+1+\frac{d+\alpha}{2};
  \frac{1}{\Vert x\Vert^2}
\right), ~~\Vert x\Vert>1
\end{array}
\right.
\end{align}
 $x\in \real^d$, where ${}_2F_1 ( a,b;c;y)$ is
 Gauss's hypergeometric function,
 see (5.2) in \cite{getoor},
Lemma~4.1 in \cite{biler2015nonlocal},
and Relation~(36) in \cite{oberman}.
\subsubsection*{Nonlinear fractional PDE}
\begin{figure}[H]
\centering
\hskip-0.4cm
\begin{subfigure}{.5\textwidth}
\centering
\includegraphics[width=\textwidth]{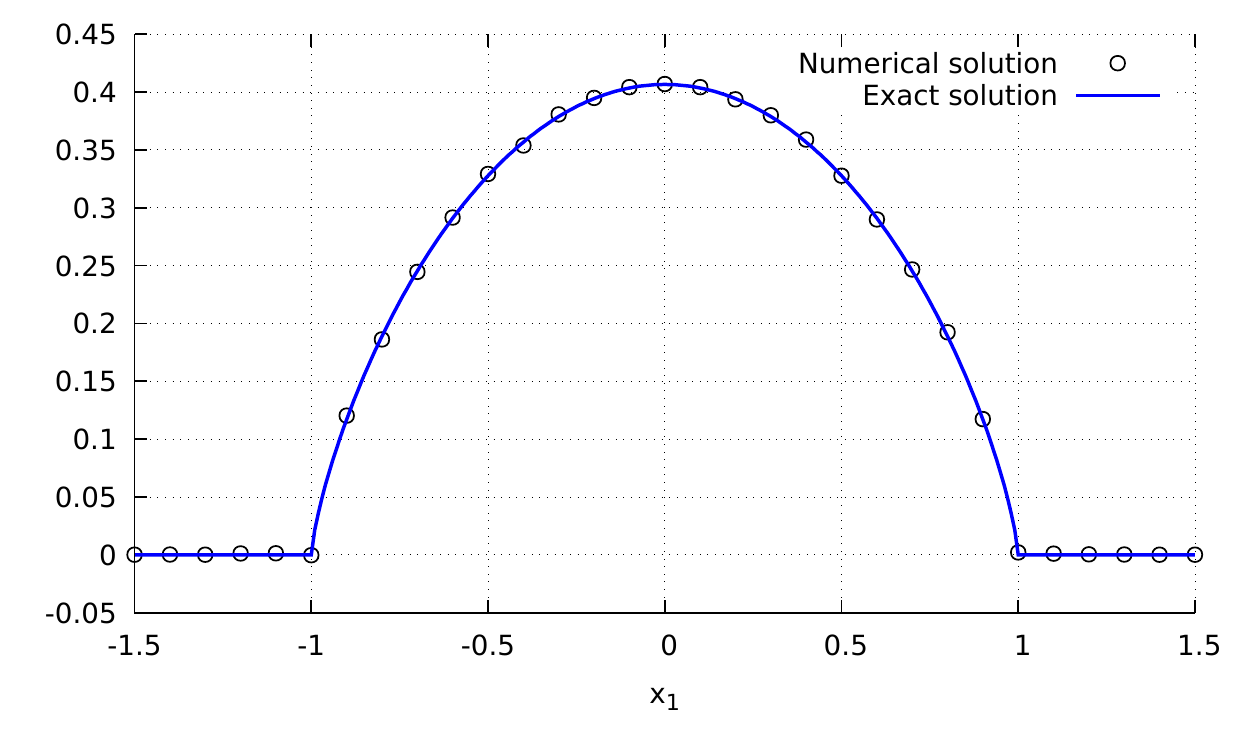}
\vskip-0.3cm
\caption{Numerical solution of \eqref{eq:nld} with $k=0$.}
\end{subfigure}
\begin{subfigure}{.50\textwidth}
\centering
\includegraphics[width=\textwidth]{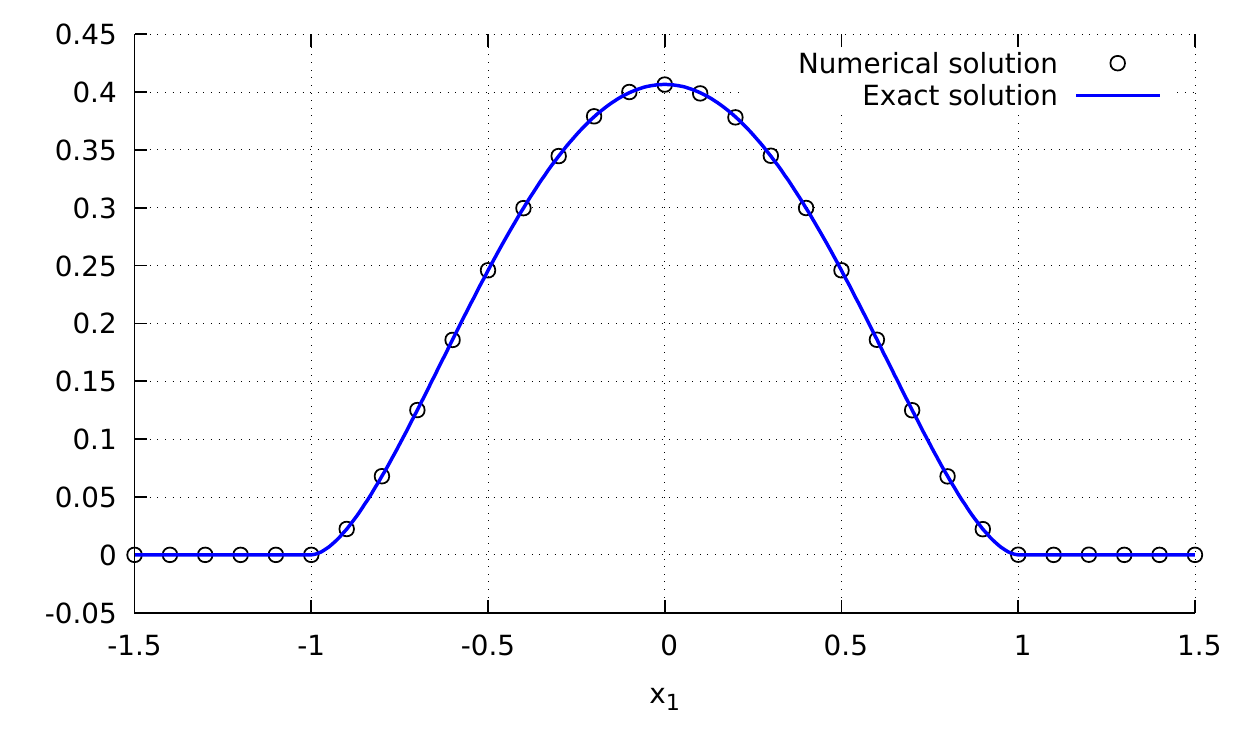}
\vskip-0.3cm
\caption{Numerical solution of \eqref{eq:nld} with $k=1$.}
\end{subfigure}
\caption{Numerical solution of \eqref{eq:nld} { in dimension $d=10$}.}
\label{fig1}
\end{figure}

Based on \eqref{pp}, we aim at recovering the explicit solution
\begin{equation}
 \label{djk1}
   u(t,x) = e^{-t}\Phi_{k,\alpha} (x)
   = e^{-t}(1-\Vert x\Vert^2)^{k+\alpha / 2}_+, \qquad
 (t,x)\in [0,T] \times \real^d.
\end{equation}
of the nonlinear PDE
\begin{equation}\label{eq:nld}
\begin{cases}
  \displaystyle
  - \frac{\partial u}{\partial t} (t,x) = \Delta_\alpha u (t,x) + e^{-t} \Psi_{k,\alpha}(x) - e^{-4t}(1-\Vert x\Vert^2)^{4k+2\alpha}_+ + u (t,x)+ u^4 (t,x),
  \medskip
  \\
u(T,x) = e^{-T}(1-\Vert x\Vert^2)^{k+\alpha / 2}_+,
\end{cases}{}
\end{equation}{}
 with $m=0$,
$f(t, x,y ) = c_0 (t,x)+ y + y^4
$
 and ${\cal L}_0 = \{ 0, 1, 4 \}$,
 $c_0 (t,x) = e^{-t} \Psi_{k,\alpha}(x) - e^{-4t}(1-\Vert x\Vert^2)^{4k+2\alpha}_+$,
 $c_1 (t,x)=c_4 (t,x)=1$. 
 The random tree associated to
 Equation~\eqref{eq:nld} started with a mark $i\in \{0,\ldots , d\}$
 branches into \textcolor{blue}{0 branch}, \textcolor{cyan}{1 branch},
 or \textcolor{violet}{4 branches}, with the mark~$0$,
 as in the following random sample:

\begin{center}
\resizebox{0.65\textwidth}{!}{
\begin{tikzpicture}[scale=0.9,grow=right, sloped]
\node[rectangle,draw,cyan,text=black,thick]{\shortstack{$t$\\$x$}}
    child {
        node[rectangle,draw,violet,text=black,thick] {\shortstack{$t+T_{\bar{1}}$\\ $X^{\bar{1}}_{T_{\bar{1}},x} $}}
            child {
                node[rectangle,draw,black,thick] {\shortstack{$T$\\ $X^{(1,4)}_{T,x}  $}}
                edge from parent
                node[above] {$(1,4)$}
                node[below]  {$0$}
            }
            child {
                node[rectangle,draw,violet,text=black,thick] {\shortstack{$t+T_{\bar{1}} + T_{(1,3)}$\\ $X^{(1,3)}_{T_{(1,3)},x} $}}
                child{
                node[rectangle,draw,blue,text=black,thick]{\shortstack{$t+T_{\bar{1}}+ T_{(1,3)}+ T_{(1,3,4)}$\\ $X^{(1,3,4)}_{T_{(1,3,4)},x} $}}
                edge from parent
                node[above]{$(1,3,4)$}
                node[below]{$0$}
                }
                child{
                node[rectangle,draw,thick]{\shortstack{$T$\\ $X^{(1,3,3)}_{T,x} $}}
                edge from parent
                node[above]{$(1,3,3)$}
                node[below]{$0$}
                }
                child{
                node[rectangle,draw,thick]{\shortstack{$T$\\ $X^{(1,3,2)}_{T,x} $}}
                edge from parent
                node[above]{$(1,3,2)$}
                node[below]{$0$}
                }
                child{
                node[rectangle,draw,thick]{\shortstack{$T$\\ $X^{(1,3,1)}_{T,x} $}}
                edge from parent
                node[above]{$(1,3,1)$}
                node[below]{$0$}
                }
                edge from parent
                node[above] {$(1,3)$}
                node[below]  {$0$}
            }
            child {
                node[rectangle,draw,blue,text=black,thick] {\shortstack{$t+T_{\bar{1}} + T_{(1,2)}$\\ $X^{(1,2)}_{T_{(1,2)},x}  $}}
                edge from parent
                node[above] {$(1,2)$}
                node[below]  {$0$}
            }
            child {
                node[rectangle,draw,cyan,text=black,thick] {\shortstack{$t+T_{\bar{1}} + T_{(1,1)}$\\ $X^{(1,1)}_{T_{(1,1)},x}  $}}
                child{
                    node[rectangle,draw,black] {\shortstack{$T$\\ $X^{(1,1,1)}_{T,x}  $}}
                    edge from parent
                    node[above]{~$(1,1,1)$~}
                    node[below]{$0$}
            }
                edge from parent
                node[above] {$(1,1)$}
                node[below]  {$0$}
            }
            edge from parent
            node[above] {$\bar{1}$}
            node[below]{$i$}
    };
\end{tikzpicture}
}
\end{center}
\noindent
 In Figure~\ref{fig1} we plot
 the numerical solutions $u(t,x_1,0,\ldots , 0)$ of \eqref{eq:nld}
 obtained from \eqref{u0tx2}
 {in terms of the first coordinate $x_1$}
 in dimension $d=10$, with $T=1$, $t=0.9$ and $\alpha =1.5$.

\subsubsection*{Nonlinear fractional PDE with gradient term}
 We consider the nonlinear PDE
 \begin{equation}
   \label{eq:gradd}
\begin{cases}
  \displaystyle
  - \frac{\partial u}{\partial t} (t,x)
  = \Delta_\alpha u(t,x)
  + u(t,x)
  + e^{-t} \Psi_{k,\alpha} (x)
  \\
  \hskip2.3cm
 \displaystyle  +
  (2k+\alpha) e^{-2t} (1-\Vert x\Vert^2)^{2k+\alpha-1}_+
  \sum_{j=1}^d x_j
  + u (t,x) \sum_{j=1}^d \frac{\partial u}{\partial x_j}(t,x)
  \medskip
  \\
u(T,x) = e^{-T}(1- \Vert x\Vert^2)^{k+\alpha / 2}_+,
\end{cases}{}
\end{equation}{}
with $m=d$,
$$
f(t, x,y , z_1,\ldots , z_d) =
c_{0,\ldots, 0}(t,x ) + y ( z_1+\cdots + z_d ),
$$
${\cal L}_d = \{ (0,\ldots, 0), (1,1,\ldots, 0), \ldots ,(1,0,\ldots, 1) \}$,
$$
  c_{0,\ldots, 0}(t,x)
 =  e^{-t} \Psi_{k,\alpha} (x) + (2k+\alpha) e^{-2t} (1-\Vert x\Vert^2)^{2k+\alpha-1}_+
( x_1+\cdots + x_d ),
$$
 and $c_{1,1,\ldots, 0}(t,x)=\cdots = c_{1,0,\ldots, 1}(t,x)=1$,
 whose explicit solution $u(t,x)$ is also given by \eqref{djk1}
according to \eqref{pp}.
In dimension $d=1$ the possible marks are $0$ and $1$,
 and the corresponding random trees
branches into
\textcolor{blue}{0 branch}, \textcolor{cyan}{1 branch},
or {2 branches}),
 as in the following random sample: 
\\

\begin{figure}[H]
\centering
\hskip-0.4cm
\begin{subfigure}{.5\textwidth}
\centering
\includegraphics[width=\textwidth]{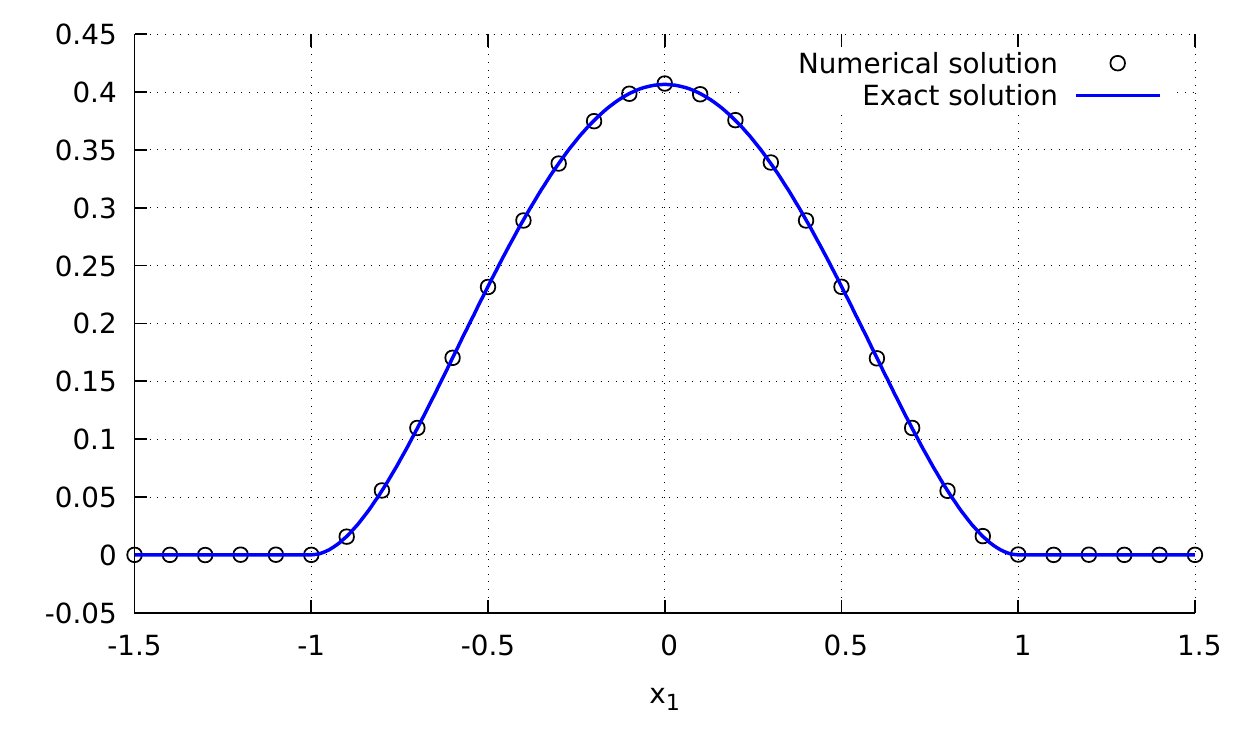}
\vskip-0.3cm
\caption{Numerical solution of \eqref{eq:gradd} with $k=1$.}
\end{subfigure}
\begin{subfigure}{.50\textwidth}
\centering
\includegraphics[width=\textwidth]{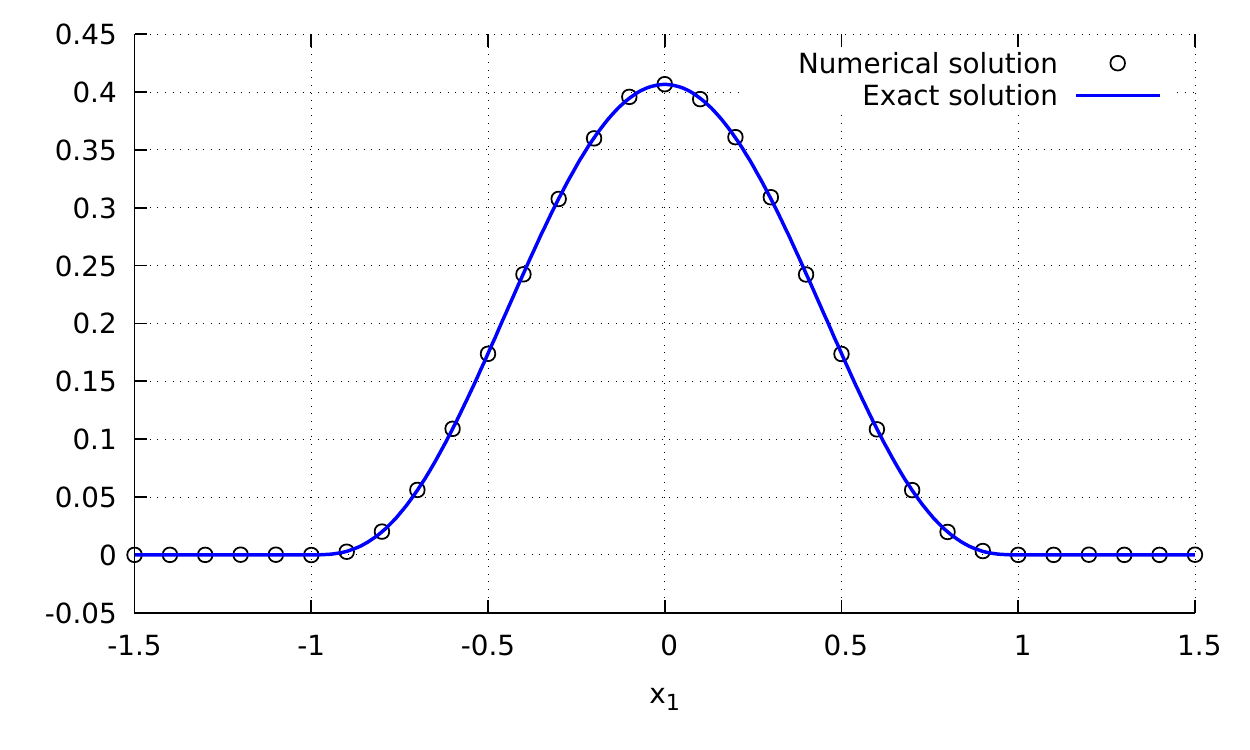}
\vskip-0.3cm
\caption{Numerical solution of \eqref{eq:gradd} with $k=2$.}
\end{subfigure}
\caption{Numerical solution of \eqref{eq:gradd} { in dimension $d=2$}.}
\label{fig2}
\end{figure}

\begin{center}
\resizebox{0.85\textwidth}{!}{
\begin{tikzpicture}[scale=1.0,grow=right, sloped]
\node[rectangle,draw,cyan,text=black,thick]{\shortstack{$t$\\$x$}}
    child {
        node[rectangle,draw,purple,text=black,thick] {\shortstack{$t+T_{\bar{1}}$\\ $X^{\bar{1}}_{T_{\bar{1}},x} $}}
            child {
                node[rectangle,draw,purple,text=black,thick] {\shortstack{$t+T_{\bar{1}} + T_{(1,2)}$\\ $X^{(1,2)}_{T_{(1,2)},x} $}}
                child{
                node[rectangle,draw,blue,text=black,thick]{\shortstack{$t+T_{\bar{1}}+ T_{(1,2)}+ T_{(1,2,2)}$\\ $X^{(1,2,2)}_{T_{(1,2,2)},x} $}}
                edge from parent
                node[above]{$(1,2,2)$}
                node[below]{$1$}
                }
                child{
                node[rectangle,draw,thick]{\shortstack{$T$\\ $X^{(1,2,1)}_{T,x} $}}
                edge from parent
                node[above]{$(1,2,1)$}
                node[below]{$0$}
                }
                edge from parent
                node[above] {$(1,2)$}
                node[below]  {$1$}
            }
            child {
                node[rectangle,draw,cyan,text=black,thick] {\shortstack{$t+T_{\bar{1}} + T_{(1,1)}$\\ $X^{(1,1)}_{T_{(1,1)},x}  $}}
                child {
                    node[rectangle,draw,blue,text=black,thick] {\shortstack{$t+T_{\bar{1}} + T_{(1,1)} + T_{(1,1,1)}$\\ $X^{(1,1,1)}_{T_{(1,1,1)},x}  $}}
                    edge from parent
                    node[above]{$(1,1,1)$}
                    node[below]{$0$}
                }
                edge from parent
                node[above] {$(1,1)$}
                node[below]  {$0$}
            }
            edge from parent
            node[above] {$\bar{1}$}
            node[below]{$0$}
    };
\end{tikzpicture}
}
\end{center}
In dimension $d>1$ the tree expands into
 $d+2$ different types of branches, namely
$0$ branch, one branch with mark $0$, and $d$ types of two branches
with one branch bearing the mark $0$ and the other branch
bearing a mark $i \in \{1,\dots,d\}$, which
 corresponds to the partial derivative with respect to $x_i$.
 In Figure~\ref{fig2} we plot
 the numerical solutions $u(t,x_1,0,\ldots , 0)$ of \eqref{eq:gradd}
 obtained from \eqref{u0tx2}
 {in terms of the first coordinate $x_1$}
 in dimension $d=2$, with $T=1$, $t=0.9$ and $\alpha =1.5$.

\subsubsection*{Fractional Burgers equation}
\noindent
 We consider the fractional Burgers equation
\begin{equation}\label{burgers_dim2}
  \frac{\partial u}{\partial t}(t,x)+ \kappa \Delta_\alpha u (t,x)- u(t,x)
  \sum_{j=1}^d
   \frac{\partial u}{\partial x_j}(t,x) = 0, \qquad
   x = (x_1,\ldots , x_d)\in \real^d,
\end{equation}
with $m=d$, $f(t, x,y , z_1,\ldots , z_d) = y ( z_1+\cdots + z_d )$
 and
 ${\cal L}_d = \{ (1,1,\ldots, 0), \ldots ,(1,0,\ldots, 1) \}$,
 $c_{1,1,\ldots, 0}(t,x)=\cdots = c_{1,0,\ldots, 1}(t,x)=1$,
 and one of the following two terminal conditions.
\begin{enumerate}[1.]
  \item Half-space terminal condition
\begin{equation}\label{cond2}
  u(T,x) = \mathbbm{1}_{[0,\infty )^d}(x_1),
    \qquad
     x = (x_1,\ldots , x_d) \in \real^d. 
\end{equation}
\item Product cosine terminal condition
\begin{equation}\label{cond1}
  u(T,x) = \cos(x_1) \cdots \cos(x_d) \mathbbm{1}_{[-\pi / 2, \pi / 2]^d}(x_1,\ldots , x_d),
  \qquad
     x = (x_1,\ldots , x_d) \in \real^d.
\end{equation}
\end{enumerate}
The random tree associated to this equation is a binary tree with $d$ types of branching. At each branching time, two branches are generated, one bearing the mark $0$ to represent $u$ and the other one bearing a mark $i \in \{1,\dots,d \} $ to represent $ \partial u / \partial x_i$,
which yields the following sample tree in dimension $d=1$.

\bigskip

\begin{center}
\resizebox{0.65\textwidth}{!}{
\begin{tikzpicture}[scale=0.9,grow=right, sloped]
\node[rectangle,draw,cyan,thick]{\shortstack{$t$\\$x$}}
    child {
        node[rectangle,draw,purple,text=black,thick] {\shortstack{$t+T_{\bar{1}}$\\ $X^{\bar{1}}_{T_{\bar{1}},x} $}}
            child {
                node[rectangle,draw,purple,text=black,thick] {\shortstack{$t+T_{\bar{1}} + T_{(1,2)}$\\ $X^{(1,2)}_{T_{(1,2)},x} $}}
                child{
                node[rectangle,draw,black,text=black,thick]{\shortstack{$T$\\ $X^{(1,2,2)}_{T,x} $}}
                edge from parent
                node[above]{$(1,2,2)$}
                node[below]{$1$}
                }
                child{
                node[rectangle,draw,thick]{\shortstack{$T$\\ $X^{(1,2,1)}_{T,x} $}}
                edge from parent
                node[above]{$(1,2,1)$}
                node[below]{$0$}
                }
                edge from parent
                node[above] {$(1,2)$}
                node[below]  {$1$}
            }
            child {
                node[rectangle,draw,black,text=black,thick] {\shortstack{$T$\\ $X^{(1,1)}_{T,x}  $}}
                edge from parent
                node[above] {$(1,1)$}
                node[below]  {$0$}
            }
            edge from parent
            node[above] {$\bar{1}$}
            node[below]{$0$}
    };
\end{tikzpicture}
}
\end{center}

\begin{figure}[H]
\centering
\hskip-0.4cm
\begin{subfigure}{.5\textwidth}
\centering
\includegraphics[width=\textwidth]{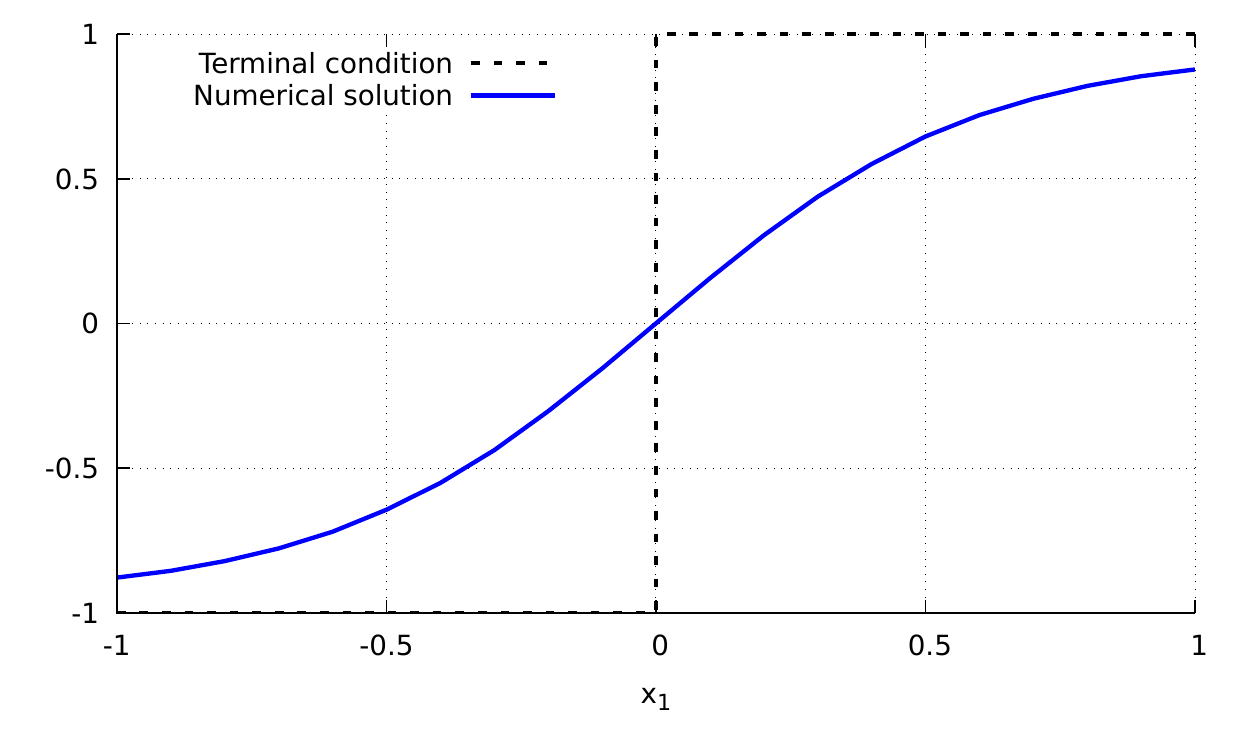}
\vskip-0.3cm
\caption{Terminal condition \eqref{cond2} and $t=0.99$.}
\end{subfigure}
\begin{subfigure}{.50\textwidth}
\centering
\includegraphics[width=\textwidth]{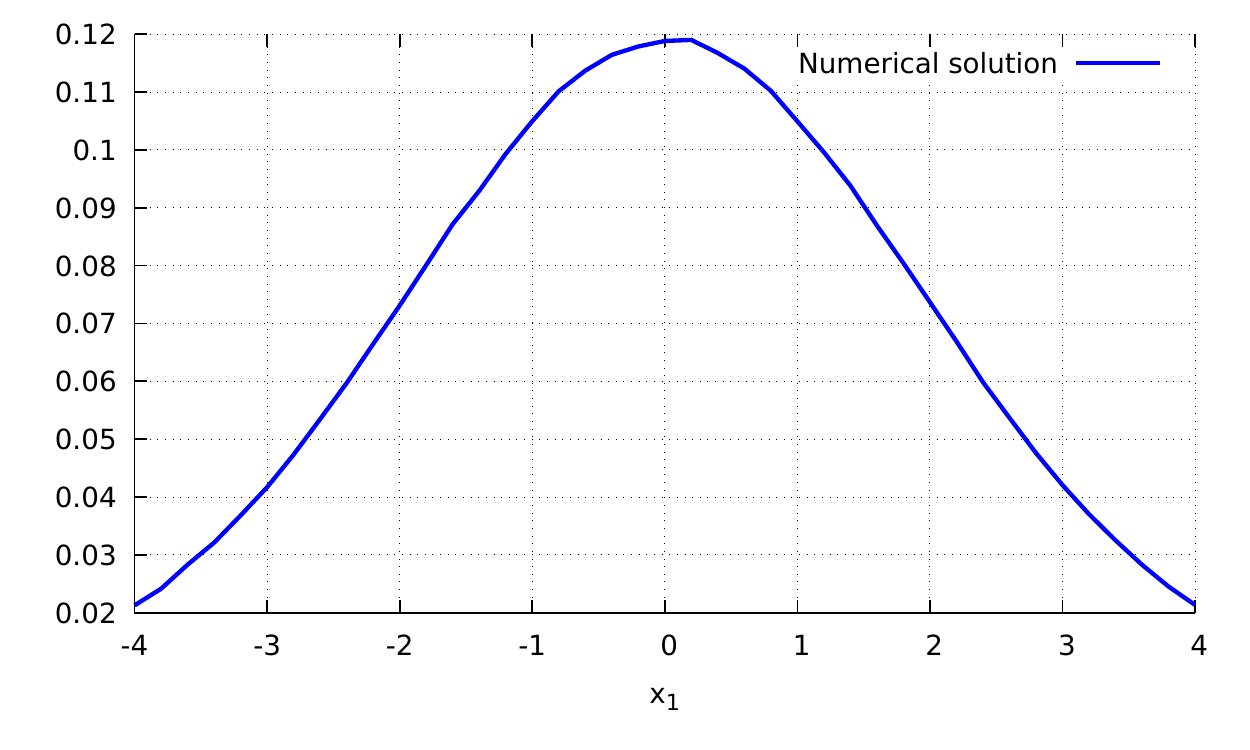}
\vskip-0.3cm
\caption{Terminal condition \eqref{cond1} and $t=0.9$.}
\end{subfigure}
\caption{Numerical solution of \eqref{burgers_dim2} { in dimension $d=2$}.}
\label{fig3}
\end{figure}

\noindent
 In Figure~\ref{fig3} we plot
 the numerical solutions $u(t,x_1,0)$ of \eqref{burgers_dim2}
 obtained from \eqref{u0tx2} 
 {in terms of the first coordinate $x_1$}
 in dimension $d=2$, with $\kappa =10$,
 $T=1$ and $\alpha =1.5$.

\subsubsection*{Data availability statement}
 No new data were created during the study. 

\footnotesize

\setcitestyle{numbers}

\def\cprime{$'$} \def\polhk#1{\setbox0=\hbox{#1}{\ooalign{\hidewidth
  \lower1.5ex\hbox{`}\hidewidth\crcr\unhbox0}}}
  \def\polhk#1{\setbox0=\hbox{#1}{\ooalign{\hidewidth
  \lower1.5ex\hbox{`}\hidewidth\crcr\unhbox0}}} \def\cprime{$'$}

\end{document}